\documentclass[a4 paper, 10pt, twosided]{article}
\usepackage{mathptmx}
\usepackage{authblk}
\usepackage{fancyhdr}
\fancypagestyle{fpg}{

\fancyhead[R]{\textit{Transactions of the Nigerian Association of Mathematical Physics}\\\textit{Volume 16, (July$-$Sept. 2021 Issue),} $pp73-98$ \\\textcopyright\textit{Trans. of NAMP}}
\fancyhead[L]{}
\fancyhead[C]{}
\fancyfoot[L]{Correspondence Author: M.E. Egwe, Email: \href{mailto:murphy.egwe@mail.ui.edu.ng}{murphy.egwe@mail.ui.edu.ng},
\textbf{\textit{Transactions of the Nigerian Association of Mathematical Physics Volume 16, (July$-$Sept. 2021), pg73$-$98}}\begin{center}\thepage\end{center}}}

\usepackage{amsmath,amssymb, amsfonts}
\usepackage{mathptmx}
\usepackage{empheq}
\usepackage{phfthm}
\usepackage{framed}
\usepackage{setspace}
\usepackage{wasysym}
\usepackage[x11names, dvipsnames]{xcolor}
\usepackage{tcolorbox}
\usepackage{pgfplots}
\usepgfplotslibrary{external}
\usepackage{tikz}
\tcbuselibrary{skins,breakable}
\usetikzlibrary{patterns}
\pgfplotsset{every axis/.append style={axis x line = middle, axis y line = middle}}
\usepackage[urlcolor=blue, colorlinks=true, linkcolor=blue, citecolor=red]{hyperref}
\usepackage{phfthm}
\usepackage{wasysym}
\makeatletter
\renewenvironment{proof}[1][\proofname]{\par
  \pushQED{\qed} \normalfont \topsep6\p@\@plus6\p@\relax
  \trivlist
  \item[\hskip\labelsep
    \normalfont\bfseries
    #1\@addpunct{.}\ignorespaces
  ]\ignorespaces
}{\popQED\endtrivlist\@endpefalse
}
\let\qed\relax 
\DeclareRobustCommand{\qed}{
  \ifmmode \mathqed
  \else
    \leavevmode\unskip\penalty\@M\hbox{}\nobreak\hspace{.5em minus .1em}
    \hbox{\qedsymbol}
  \fi
}
\makeatother
\allowdisplaybreaks
\title{\textbf{New Integral Representations for the Bessel Functions}}
\author[1]{Abdulhafeez A. Abdulsalam}
\author[2]{M.E. Egwe}
\affil[1,2]{Department of Mathematics, University of Ibadan, Ibadan, Nigeria} 
\affil[ ]{$^{1}$\href{mailto:aabdulsalam030@stu.ui.edu.ng}{aabdulsalam030@stu.ui.edu.ng}, \href{mailto:hafeez147258369@gmail.com}{hafeez147258369@gmail.com}, $^{2}$\href{mailto:murphy.egwe@mail.ui.edu.ng}{murphy.egwe@mail.ui.edu.ng}}

\newcommand{\rmd}{\mathrm}
\newcommand{\om}{\omega}
\newcommand{\var}{\vartheta}
\DeclareMathAlphabet{\mathcal}{OMS}{cmsy}{m}{n}
\definecolor{lightgray}{gray}{0.75}
\date{}
\begin{document}

\maketitle
\thispagestyle{fpg}

\vspace{-3mm} \textcolor{blue}{\hrulefill}\\
\vspace{-3mm}
\small
\begin{abstract}
In this paper, new integral representations for the Bessel $J$ and $I$ functions were presented and their results were used to derive an expression for the Modified Bessel $K$ function. 
\end{abstract}
\vspace{-3mm} \textcolor{blue}{\hrulefill}\\
\vspace{-3mm}

\section{Introduction}

The theory of Bessel functions is intimately connected with the theory of a certain type of differential equation of the first order, known as Riccati's equation. In fact, a Bessel function is usually defined as a particular solution of a linear differential equation of the second order (known as Bessel's equation) which is derived from Riccati's equation by elementary transformation. \cite{Treat2}\\

Bessel functions, first defined by the mathematician Daniel Bernoulli and then generalized by Friedrich Bessel, are canonical solutions $y(x)$ of Bessel's differential equation

\begin{center}
\tcbox[interior style={opacity=0.0},extrude left by=-0.2cm]{$\displaystyle x^2\frac{\rmd{d}^2y}{\rmd{d}x^2} + x\frac{\rmd{d}y}{\rmd{d}x} + (x^2 - \alpha^2)y = 0$}
\end{center}

Bessel functions of the first kind, denoted as $J_{\alpha}(x)$, are solutions of Bessel's differential equation. For integer or positive $\alpha$, Bessel functions of the first kind are finite at the origin; while for negative non-integer $\alpha$, Bessel functions of the first kind diverge as $x$ approaches zero. It is possible to define the function by its series expansion around $x=0$, which can be found by applying the Frobenius method to Bessel's equation 

$$J_{\alpha}(x) = \sum_{m=0}^{\infty} \frac{(-1)^m}{m!\Gamma(m + \alpha + 1)} \left(\frac{x}{2}\right)^{2m + \alpha}\hspace{0.5cm}\forall x \in \mathbb{C}$$

where $\Gamma(z)$ is the gamma function, a shifted generalization of the factorial function to non-integer values. \cite{Treat5}

\section{Integral representations for $J_{\alpha}(x)$ and $J_{0}(x)$}

$$J_{\alpha}(x) = \sum_{m=0}^{\infty} \frac{(-1)^m}{m!\Gamma(m + \alpha + 1)} \left(\frac{x}{2}\right)^{2m + \alpha}$$

\begin{align*}
J_{\alpha}(x) &= \sum_{m=0}^{\infty} \frac{(-1)^m}{m!\Gamma(m + \alpha + 1)} \left(\frac{x}{2}\right)^{2m + \alpha}
\\&= \frac{1}{\Gamma(\alpha)} \sum_{m=0}^{\infty} \frac{(-1)^m\Gamma(m + 1)\Gamma(\alpha)}{m!\Gamma(m + 1)\Gamma(m + \alpha + 1)} \left(\frac{x}{2}\right)^{2m + \alpha} 
\\&= \frac{1}{\Gamma(\alpha)}\sum_{m=0}^{\infty} \frac{(-1)^m\Gamma(m + 1)\Gamma(\alpha)}{\left(m!\right)^2\Gamma(m + \alpha + 1)} \left(\frac{x}{2}\right)^{2m + \alpha}
\\&= \frac{1}{\Gamma(\alpha)}\sum_{m=0}^{\infty} \frac{(-1)^m}{\left(m!\right)^2}\beta(m + 1, \alpha) \left(\frac{x}{2}\right)^{2m + \alpha}
\\&= \frac{1}{\Gamma(\alpha)}\sum_{m=0}^{\infty} \frac{(-1)^m}{\left(m!\right)^2} \int_0^1 t^{m} \left(1 - t\right)^{\alpha - 1} \rmd{d}t \left(\frac{x}{2}\right)^{2m + \alpha}
\\&= \frac{1}{\Gamma(\alpha)}\left(\frac{x}{2}\right)^{\alpha} \int_0^1 \sum_{m=0}^{\infty} \frac{(-1)^m}{\left(m!\right)^2} t^{m}  \left(\frac{x}{2}\right)^{2m} \left(1 - t\right)^{\alpha - 1} \rmd{d}t
\\&= \frac{1}{\Gamma(\alpha)}\left(\frac{x}{2}\right)^{\alpha} \int_0^1 \sum_{m=0}^{\infty} \frac{(-1)^m}{\left(m!\right)^2}  \left(\frac{x\sqrt{t}}{2}\right)^{2m} \left(1 - t\right)^{\alpha - 1} \rmd{d}t
\end{align*}

We notice that 
\begin{align*}
J_{0}(x) &= \sum_{m=0}^{\infty} \frac{(-1)^m}{m!\Gamma(m + 0 + 1)} \left(\frac{x}{2}\right)^{2m +0} 
\\&= \sum_{m=0}^{\infty} \frac{(-1)^m}{m!\Gamma(m + 1)} \left(\frac{x}{2}\right)^{2m} 
\\&= \sum_{m=0}^{\infty} \frac{(-1)^m}{\left(m!\right)^2} \left(\frac{x}{2}\right)^{2m} 
\end{align*}

\begin{align}
\implies J_{\alpha}(x) &= \frac{1}{\Gamma(\alpha)}\left(\frac{x}{2}\right)^{\alpha} \int_0^1 J_{0}(x\sqrt{t})\left(1 - t\right)^{\alpha - 1} \rmd{d}t\nonumber
\\&= \frac{1}{\Gamma(\alpha)}\left(\frac{x}{2}\right)^{\alpha} \int_0^1 J_{0}\left(x\sqrt{1 - t}\right)t^{\alpha - 1} \rmd{d}t\label{refe}
\end{align}

By Legendre's duplication formula,

$$\Gamma\left(m + \frac{1}{2}\right)= \frac{\Gamma(2m)\sqrt{\pi}}{\Gamma(m)2^{2m - 1}}$$

$$\implies \frac{\Gamma\left(m + \frac{1}{2}\right)2^{2m - 1}}{\Gamma(2m)\sqrt{\pi}} = \frac{1}{\Gamma(m)}$$

\begin{align*}
\implies J_{0}(x) &=  \sum_{m=0}^{\infty} \frac{(-1)^m}{m!\Gamma(m + 1)} \left(\frac{x}{2}\right)^{2m}
\\&= \sum_{m=0}^{\infty} \frac{(-1)^m}{\Gamma(m + 1)^2} \left(\frac{x}{2}\right)^{2m}
\\&= \sum_{m=1}^{\infty} \frac{(-1)^m}{\Gamma(m)\Gamma(m)} \left(\frac{x}{2}\right)^{2m-2}
\\&= \sum_{m = 1}^{\infty} \frac{(-1)^{m-1}\Gamma\left(m + \frac{1}{2}\right)2^{2m - 1}}{\Gamma(2m)\Gamma(m)\sqrt{\pi}} \left(\frac{x}{2}\right)^{2m - 2}
\\&= \sum_{m = 1}^{\infty} \frac{2m(-1)^{m-1}\,\Gamma\left(m + \frac{1}{2}\right)2^{2m - 1}}{2m\Gamma(2m)\Gamma(m)\sqrt{\pi}} \left(\frac{x}{2}\right)^{2m - 2} = \sum_{m = 1}^{\infty} \frac{2m(-1)^{m-1}\,\Gamma\left(m + \frac{1}{2}\right)2^{2m - 1}}{\Gamma(2m + 1)\Gamma(m)\sqrt{\pi}} \left(\frac{x}{2}\right)^{2m - 2}
\\&=\sum_{m = 1}^{\infty} \frac{2m(-1)^{m-1}\,\Gamma\left(m + \frac{1}{2}\right)\Gamma\left(m + \frac{1}{2}\right)}{\Gamma(2m + 1)} \cdot \frac{2^{2m - 1}}{\Gamma\left(m + \frac{1}{2}\right)\Gamma(m)\sqrt{\pi}}\left(\frac{x}{2}\right)^{2m - 2}
\\&=\sum_{m = 1}^{\infty} 2m(-1)^{m-1}\,\beta\left(m + \frac{1}{2}, m + \frac{1}{2}\right) \cdot \frac{2^{2m - 1}}{\Gamma\left(m + \frac{1}{2}\right)\Gamma(m)\sqrt{\pi}}\left(\frac{x}{2}\right)^{2m - 2}
\\&=\sum_{m = 1}^{\infty} 2m(-1)^{m-1}\,\beta\left(m + \frac{1}{2}, m + \frac{1}{2}\right) \cdot \frac{2^{2m - 1}}{\frac{\Gamma(2m)\sqrt{\pi}}{\Gamma(m)2^{2m - 1}}\cdot\Gamma(m)\sqrt{\pi}}\left(\frac{x}{2}\right)^{2m - 2}
\\&=\sum_{m = 1}^{\infty} 2m(-1)^{m-1}\,\beta\left(m + \frac{1}{2}, m + \frac{1}{2}\right) \cdot \frac{2^{4m - 2}}{\pi\Gamma(2m)}\left(\frac{x}{2}\right)^{2m - 2}
\\&= \frac{1}{2\pi}\sum_{m = 1}^{\infty} m(-1)^{m-1}\,\beta\left(m + \frac{1}{2}, m + \frac{1}{2}\right) \cdot \frac{1}{\Gamma(2m)}\left(2x\right)^{2m}\left(\frac{x}{2}\right)^{-2}
\\&= \frac{1}{2\pi}\sum_{m = 1}^{\infty} m(-1)^{m-1}\int_0^1 \om^{m - \frac{1}{2}}\left(1 - \om\right)^{m - \frac{1}{2}} \rmd{d}\om \cdot \frac{1}{\Gamma(2m)}\left(2x\right)^{2m}\left(\frac{x}{2}\right)^{-2} 
\\&= \frac{2}{\pi x^2}\int_0^1\left(\om^{-\frac{1}{2}}\left(1 - \om\right)^{-\frac{1}{2}} \sum_{m = 1}^{\infty} m(-1)^{m-1}\om^{m}\left(1 - \om\right)^{m} \cdot \frac{1}{\Gamma(2m)}\left(2x\right)^{2m}\right)\rmd{d}\om
\\&= \frac{2}{\pi x^2}\int_0^1\left(\om^{-\frac{1}{2}}\left(1 - \om\right)^{-\frac{1}{2}} \sum_{m = 1}^{\infty}\frac{(-1)^{m-1} m \left(4x^2\left(\om - \om^2\right)\right)^m}{\Gamma(2m)}\right)\rmd{d}\om
\end{align*}

Consider 
$$\sum_{m = 1}^{\infty}\frac{(-1)^{m-1} m z^m}{\Gamma(2m)}$$

\begin{align*}
\sum_{m = 1}^{\infty}\frac{(-1)^{m-1} m z^m}{\Gamma(2m)} &= z\sum_{m = 1}^{\infty}\frac{(-1)^{m-1}\frac{\partial}{\partial z}\left(z^m\right)}{\Gamma(2m)} = z\frac{\partial}{\partial z}\left(\sum_{m = 1}^{\infty}\frac{(-1)^{m-1} z^m}{\Gamma(2m)}\right) =  z\frac{\partial}{\partial z}\left(\sum_{m = 0}^{\infty}\frac{(-1)^{m} z^{m+1}}{\Gamma(2m + 2)}\right)
\\&= z\frac{\partial}{\partial z}\left(\sum_{m = 0}^{\infty}\frac{(-1)^{m} z^{m+1}}{(2m + 1)!}\right) =  z\frac{\partial}{\partial z}\left(\sqrt{z}\sum_{m = 0}^{\infty}\frac{(-1)^{m} \left(\sqrt{z}\right)^{2m+1}}{(2m + 1)!}\right) 
\\&= z\frac{\partial}{\partial z}\left(\sqrt{z}\sin{\sqrt{z}}\right) = z\left(\frac{1}{2\sqrt{z}} \sin{\sqrt{z}} + \frac{1}{2}\cos{\sqrt{z}}\right) = \frac{1}{2}\left(\sqrt{z}\sin{\sqrt{z}} + z\cos{\sqrt{z}}\right)
\end{align*}

\begin{align*}
\implies J_{0}(x) &=   \frac{2}{\pi x^2}\int_0^1\left(\om^{-\frac{1}{2}}\left(1 - \om\right)^{-\frac{1}{2}} \sum_{m = 1}^{\infty}\frac{(-1)^m m \left(4x^2\left(\om - \om^2\right)\right)^m}{\Gamma(2m)}\right)\rmd{d}\om
\\&= \frac{2}{\pi x^2}\cdot \frac{1}{2} \int_0^1 \left(\om^{-\frac{1}{2}}\left(1 - \om\right)^{-\frac{1}{2}}\left(\sqrt{4x^2\om(1 - \om)}\sin\left(\sqrt{4x^2\om(1 - \om)}\right)\right.\right.
\\&\hspace{3cm}\left.\left. + 4x^2\om(1 - \om)\cos\left(\sqrt{4x^2\om(1 - \om)}\right)\right)\right)\rmd{d}\om
\\&= \frac{1}{\pi x^2}\int_0^1 2x\sin\left(\sqrt{4x^2\om(1 - \om)}\right) +  4x^2\sqrt{\om(1 - \om)}\cos\left(\sqrt{4x^2\om(1 - \om)}\right) \rmd{d}\om
\end{align*}

On the transformation $\om \rightarrow \sin^2\theta$,
$$J_{0}(x) = \frac{1}{\pi x^2}\int_0^{\frac{\pi}{2}} \left(2x\sin\left(2x\sin\theta\cos\theta\right) + 4x^2\sin\theta\cos\theta\cos\left(2x\sin\theta\cos\theta\right)\right) 2\sin\theta\cos\theta \rmd{d}\theta$$

\begin{align*}
\implies J_{0}(x) &= \frac{1}{\pi x^2}\int_0^{\frac{\pi}{2}}\left(2x\sin\left(2x\sin\theta\cos\theta\right) + 4x^2\sin\theta\cos\theta\cos\left(2x\sin\theta\cos\theta\right)\right)2\sin\theta\cos\theta\,\rmd{d}\theta
\\&=\frac{1}{\pi x^2}\int_0^{\frac{\pi}{2}}\left(2x\sin\left(x\sin{2\theta}\right) + 2x^2\sin{2\theta}\cos\left(x\sin{2\theta}\right)\right)\sin{2\theta}\,\rmd{d}\theta
\\&=\frac{1}{\pi x^2}\int_0^{\pi}\left(x\sin\left(x\sin{\phi}\right) + x^2\sin{\phi}\cos\left(x\sin{\phi}\right)\right)\sin{\phi}\,\rmd{d}\phi
\\&=\frac{1}{\pi x}\int_0^{\pi}\left(\sin\left(x\sin{\phi}\right) + x\sin{\phi}\cos\left(x\sin{\phi}\right)\right)\sin{\phi}\,\rmd{d}\phi
\end{align*}

Therefore, 
\begin{align}
J_{0}(x) &=\sum_{m=0}^{\infty} \frac{(-1)^m}{m!\Gamma(m + 1)} \left(\frac{x}{2}\right)^{2m}
\\&=\frac{1}{\pi x}\int_0^{\pi}\left(\sin\left(x\sin{\phi}\right) + x\sin{\phi}\cos\left(x\sin{\phi}\right)\right)\sin{\phi}\,\rmd{d}\phi\label{imag}
\end{align}

\begin{align*}
J_{\alpha}(x)  &= \frac{1}{\Gamma(\alpha)}\left(\frac{x}{2}\right)^{\alpha} \int_0^1 J_{0}\left(x\sqrt{1 - t}\right)t^{\alpha - 1} \rmd{d}t
\\&= \frac{1}{\Gamma(\alpha)}\left(\frac{x}{2}\right)^{\alpha} \int_0^{\frac{\pi}{2}} 2\sin^{2\alpha - 1}{\var}\cos{\var}J_{0}\left(x\cos{\var}\right) \rmd{d}\var
\\&= \frac{1}{\pi x\Gamma(\alpha)}\left(\frac{x}{2}\right)^{\alpha} \int_0^{\frac{\pi}{2}} \int_0^{\pi} \frac{1}{\cos{\var}}\biggr[2\sin^{2\alpha - 1}{\var}\sin{\phi}\cos{\var} \left[\sin\left(x\cos\var\sin{\phi}\right) 
\right.
\\&\left.\hspace{7cm} + x\cos\var\sin{\phi}\cos\left(x\cos\var\sin{\phi}\right)\right]\biggr]\,\rmd{d}\phi \rmd{d}\var
\\&= \frac{1}{\pi x\Gamma(\alpha)}\left(\frac{x}{2}\right)^{\alpha} \int_0^{\frac{\pi}{2}} \int_0^{\pi}\biggr(2\sin^{2\alpha - 1}{\var}\sin{\phi}\left[\sin\left(x\cos\var\sin{\phi}\right) + x\cos\var\sin{\phi}\cos\left(x\cos\var\sin{\phi}\right)\right]\biggr)\,\rmd{d}\phi \rmd{d}\var
\end{align*}

\section{Using the new integral representation to establish Lipschitz's result}

It was shown by Lipschitz that \cite{Treat3}
$$\int_0^\infty e^{-at} J_0(bt) \,\, \rmd{d}t = \frac{1}{\sqrt{a^2 + b^2}}$$

The simplest method of establishing this result is to replace the Bessel coefficient by Parseval's integral representation and then change the order of integration by Fubini's theorem.\\

However, implementing the integral representation I discovered, Lipschitz's result can also be established.

\begin{align*}
J_{0}(x) &= \frac{1}{\pi x}\int_0^{\pi}\left(\sin\left(x\sin{\phi}\right) + x\sin{\phi}\cos\left(x\sin{\phi}\right)\right)\sin{\phi}\,\rmd{d}\phi
\\\implies J_{0}(bt) &= \frac{1}{\pi bt}\int_0^{\pi}\left(\sin\left(bt\sin{\phi}\right) + bt\sin{\phi}\cos\left(bt\sin{\phi}\right)\right)\sin{\phi}\,\rmd{d}\phi
\end{align*}

\begin{align*}
\int_0^\infty e^{-at} J_0(bt) \,\, \rmd{d}t  &=  \int_0^\infty\frac{e^{-at}}{\pi bt}\int_0^{\pi}\left(\sin\left(bt\sin{\phi}\right) + bt\sin{\phi}\cos\left(bt\sin{\phi}\right)\right)\sin{\phi}\,\rmd{d}\phi\,\rmd{d}t
\\&= \int_0^{\pi} \int_0^\infty e^{-at}\left(\frac{\sin\left(bt\sin{\phi}\right)}{\pi bt} + \frac{1}{\pi}\sin{\phi}\cos\left(bt\sin{\phi}\right)\right)\sin{\phi}\,\rmd{d}t\,\rmd{d}\phi
\\&= \int_0^{\pi}\left(\int_0^\infty e^{-at}\frac{\sin\left(bt\sin{\phi}\right)\sin{\phi}}{\pi bt}\rmd{d}t  + \int_0^\infty\frac{e^{-at}}{\pi}\sin^2{\phi}\cos\left(bt\sin{\phi}\right)\,\rmd{d}t\right)\rmd{d}\phi
\end{align*}

Let
$$\Phi(b,\phi) = \int_0^\infty e^{-at}\frac{\sin\left(bt\sin{\phi}\right)}{\pi bt}\rmd{d}t$$

\begin{align*}
\Phi(b,\phi) &= \frac{1}{\pi b} \int_0^\infty e^{-at}\frac{\sin\left(bt\sin{\phi}\right)}{t}\rmd{d}t
\\&= \frac{1}{\pi b} \int_0^\infty e^{-at}\int_0^{b\sin{\phi}}\cos(\tau t)\, \rmd{d}\tau\, \rmd{d}t
\\&= \frac{1}{\pi b}\int_0^{b\sin{\phi}} \int_0^\infty e^{-at}\cos(\tau t) \,\rmd{d}t\, \rmd{d}\tau
\\&= \frac{1}{\pi b}\int_0^{b\sin{\phi}} \mathcal{L}\left\{\cos(\tau t)\right\}\left(a\right)\rmd{d}\tau
\\&= \frac{1}{\pi b}\int_0^{b\sin{\phi}} \frac{a}{\tau^2 + a^2} \rmd{d}\tau
\\&= \frac{1}{\pi b} \arctan\left(\frac{\tau}{a}\right) \biggr\vert_0^{b\sin{\phi}}
\\&= \frac{1}{\pi b} \arctan\left(\frac{b\sin{\phi}}{a}\right) 
\end{align*}

Let 
$$\eta(b,\phi) = \int_0^\infty\frac{e^{-at}}{\pi}\cos\left(bt\sin{\phi}\right)\,\rmd{d}t$$

\begin{align*}
\eta(b,\phi) &= \int_0^\infty\frac{e^{-at}}{\pi}\cos\left(bt\sin{\phi}\right)\,\rmd{d}t
\\&= \frac{1}{\pi} \int_0^\infty e^{-at}\cos\left(bt\sin{\phi}\right)\,\rmd{d}t
\\&= \frac{1}{\pi} \mathcal{L}\left\{\cos\left(bt\sin{\phi}\right)\right\}\left(a\right)
\\&= \frac{1}{\pi} \frac{a}{b^2\sin^2{\phi} + a^2}
\end{align*}

\begin{align*}
\int_0^\infty e^{-at} J_0(bt) \,\, \rmd{d}t &= \int_0^{\pi}\left(\Phi(b, \phi)\sin{\phi}  + \eta(b, \phi)\sin^2{\phi}\right)\rmd{d}\phi
\\&= \frac{1}{\pi b} \int_0^\pi \arctan\left(\frac{b\sin{\phi}}{a}\right)  \sin{\phi}\, \rmd{d}\phi + \frac{1}{\pi}\int_0^\pi \frac{a\sin^2{\phi}}{b^2\sin^2{\phi} + a^2} \rmd{d}\phi
\end{align*}

\begin{align*}
\frac{1}{\pi}\int_0^\pi \frac{a\sin^2{\phi}}{b^2\sin^2{\phi} + a^2} \rmd{d}\phi &= \frac{a}{\pi b^2}\int_0^\pi \frac{b^2\sin^2{\phi}}{b^2\sin^2{\phi} + a^2} \rmd{d}\phi 
\\&= \frac{a}{\pi b^2}\int_0^\pi \frac{b^2\sin^2{\phi} + a^2 - a^2}{b^2\sin^2{\phi} + a^2} \rmd{d}\phi 
\\&= \frac{a}{\pi b^2}\int_0^\pi \rmd{d}\phi  - \frac{a}{\pi b^2}\int_0^\pi \frac{a^2}{b^2\sin^2{\phi} + a^2} \rmd{d}\phi 
\\&= \frac{\pi a}{\pi b^2}  - \frac{a^3}{\pi b^2}\int_0^\pi \frac{1}{b^2\sin^2{\phi} + a^2} \rmd{d}\phi 
\end{align*}

\begin{align*}
\int_0^\pi \frac{1}{b^2\sin^2{\phi} + a^2} \rmd{d}\phi  &{\stackrel{\nu \rightarrow \tan\left(\frac{\phi}{2}\right)}{=}}\int_0^\infty\frac{1}{b^2\left(\frac{2\nu}{1 + \nu^2}\right)^2 + a^2} \frac{2\rmd{d}\nu}{1 + \nu^2} 
\\&= 2\int_0^\infty\frac{1 + \nu^2}{4b^2\nu^2 + a^2\left(1+ \nu^2\right)^2} \rmd{d}\nu
\\&= 2\int_0^\infty\left(\frac{1}{2a(a\nu^2 + a + 2ib\nu)}+ \frac{1}{2a(a\nu^2 + a - 2ib\nu)}\right) \rmd{d}\nu
\\&= \frac{1}{a^2}\int_0^\infty\left(\frac{1}{\nu^2 + 1 + \frac{2ib\nu}{a}}+ \frac{1}{\nu^2 + 1 - \frac{2ib\nu}{a}}\right) \rmd{d}\nu
\\&= \frac{1}{a^2}\int_0^\infty\left(\frac{1}{\left(\nu + \frac{ib}{a}\right)^2 +  1 - \left(\frac{ib}{a}\right)^2} + \frac{1}{\left(\nu - \frac{ib}{a}\right)^2 + 1 - \left(\frac{ib}{a}\right)^2}\right) \rmd{d}\nu
\\&= \frac{1}{a^2}\int_0^\infty\left(\frac{1}{\left(\nu + \frac{ib}{a}\right)^2 + 1 + \left(\frac{b}{a}\right)^2}+\frac{1}{\left(\nu - \frac{ib}{a}\right)^2 + 1 + \left(\frac{b}{a}\right)^2}\right) \rmd{d}\nu
\\&= \frac{1}{a^2}\cdot\frac{1}{\sqrt{1 + \frac{b^2}{a^2}}}\left(\arctan\left(\frac{\nu + \frac{ib}{a}}{\sqrt{1 + \frac{b^2}{a^2}}}\right) + \arctan\left(\frac{\nu - \frac{ib}{a}}{\sqrt{1 + \frac{b^2}{a^2}}}\right)\right)\biggr\vert_0^{\infty}
\\&=  \frac{1}{a^2}\cdot\frac{1}{\sqrt{1 + \frac{b^2}{a^2}}}\left(\frac{\pi}{2} + \frac{\pi}{2}\right)
\\&= \frac{\pi}{a^2}\cdot\frac{1}{\sqrt{1 + \frac{b^2}{a^2}}} = \frac{\pi}{a\sqrt{a^2 + b^2}}
\end{align*}

\begin{align*}
\therefore \frac{1}{\pi}\int_0^\pi \frac{a\sin^2{\phi}}{b^2\sin^2{\phi} + a^2} \rmd{d}\phi &= \frac{\pi a}{\pi b^2}  - \frac{a^3\pi}{a\pi b^2\sqrt{a^2 + b^2}}
\\&= \frac{a}{b^2}  - \frac{a^2}{b^2\sqrt{a^2 + b^2}}
\end{align*}\leavevmode\newline

\begin{align*}
\frac{1}{\pi b} \int_0^\pi \arctan\left(\frac{b\sin{\phi}}{a}\right)  \sin{\phi}\, \rmd{d}\phi &= \frac{1}{\pi b} \int_0^\pi \arctan\left(\frac{b\sin{\phi}}{a}\right)  \rmd{d}\left(-\cos{\phi}\right)
\\&= -\frac{\cos{\phi}}{\pi b}\arctan\left(\frac{b\sin{\phi}}{a}\right)\biggr\vert_0^\pi + \frac{1}{\pi b} \int_0^\pi \frac{\frac{b\cos{\phi}}{a}}{1 + \left(\frac{b\sin{\phi}}{a}\right)^2}\cos{\phi}\rmd{d}\phi
\\&= \frac{1}{\pi a} \int_0^\pi \frac{\cos^2{\phi}}{1 + \left(\frac{b\sin{\phi}}{a}\right)^2}\rmd{d}\phi
\\&= \frac{a}{\pi} \int_0^\pi \frac{\cos^2{\phi}}{a^2 + b^2\sin^2{\phi}}\rmd{d}\phi
\end{align*}

But
$$\int_0^\pi \frac{1}{b^2\sin^2{\phi} + a^2} \rmd{d}\phi = \int_0^\infty\frac{1}{b^2\left(\frac{2\nu}{1 + \nu^2}\right)^2 + a^2} \frac{2\rmd{d}\nu}{1 + \nu^2}  = \frac{\pi}{a\sqrt{a^2 + b^2}}$$

\begin{align*}
\implies \frac{a}{\pi} \int_0^\pi \frac{\cos^2{\phi}}{a^2 + b^2\sin^2{\phi}}\rmd{d}\phi &= \frac{a}{\pi} \int_0^\pi \frac{1 - \sin^2{\phi}}{a^2 + b^2\sin^2{\phi}}\rmd{d}\phi
\\&= \frac{a}{\pi}\int_0^\pi \frac{1}{b^2\sin^2{\phi} + a^2} \rmd{d}\phi  - \frac{1}{\pi} \int_0^\pi \frac{a\sin^2{\phi}}{a^2 + b^2\sin^2{\phi}}\rmd{d}\phi
\\&= \frac{a}{\pi}\cdot\frac{\pi}{a\sqrt{a^2 + b^2}} - \left(\frac{a}{b^2}  -  \frac{a^2}{b^2\sqrt{a^2 + b^2}}\right)
\\&= \frac{1}{\sqrt{a^2 + b^2}} - \frac{a}{b^2}  + \frac{a^2}{b^2\sqrt{a^2 + b^2}} 
\end{align*}

\begin{align*}
\implies \int_0^\infty e^{-at} J_0(bt) \,\, \rmd{d}t &= \int_0^{\pi}\left(\Phi(b, \phi)\sin{\phi}  + \eta(b, \phi)\sin^2{\phi}\right)\rmd{d}\phi
\\&= \frac{1}{\pi b} \int_0^\pi \arctan\left(\frac{b\sin{\phi}}{a}\right)  \sin{\phi}\, \rmd{d}\phi + \frac{1}{\pi}\int_0^\pi \frac{a\sin^2{\phi}}{b^2\sin^2{\phi} + a^2} \rmd{d}\phi
\\&= \frac{1}{\sqrt{a^2 + b^2}} - \frac{a}{b^2}  + \frac{a^2}{b^2\sqrt{a^2 + b^2}} + \frac{a}{b^2}  - \frac{a^2}{b^2\sqrt{a^2 + b^2}}
\\&= \frac{1}{\sqrt{a^2 + b^2}}
\end{align*}

Hence, 
$$\int_0^\infty e^{-at} J_0(bt) \,\, \rmd{d}t = \frac{1}{\sqrt{a^2 + b^2}}$$

\section{ \normalfont \bf Evaluating $\displaystyle\int_0^{\infty} e^{-at}J_{\alpha}(bt)t^{-\alpha}\,\rmd{d}t$,  $\displaystyle\int_0^{\infty}J_{\frac{1}{2}}(bt)t^{-\frac{1}{2}}\,\rmd{d}t$, \\ $\displaystyle\int_0^{\infty}e^{-at} J_{\frac{1}{2}}(bt)t^{-\frac{1}{2}}\,\rmd{d}t$, $\displaystyle\int_0^{\infty}e^{-at} J_{1}(bt)t^{-1}\,\rmd{d}t$, $\displaystyle\int_0^{\infty} J_{1}(bt)t^{-1}\,\rmd{d}t$ and $\displaystyle\int_0^{\infty}J_{\alpha}(bt)t^{-\alpha}\,\rmd{d}t$ with the new integral representation}

\subsection{$\displaystyle\int_0^{\infty} e^{-at}J_{\alpha}(bt)t^{-\alpha}\,\rmd{d}t$}
It has been established previously that
\begin{align*}
J_{\alpha}(x) &= \frac{1}{\pi x\Gamma(\alpha)}\left(\frac{x}{2}\right)^{\alpha} \int_0^{\frac{\pi}{2}} \int_0^{\pi}\biggr(2\sin^{2\alpha - 1}{\var}\sin{\phi}\left[\sin\left(x\cos\var\sin{\phi}\right) \right.
\\&\hspace{7cm}\left. + x\cos\var\sin{\phi}\cos\left(x\cos\var\sin{\phi}\right)\right]\biggr)\,\rmd{d}\phi \rmd{d}\var
\end{align*}

\begin{align*}
\implies J_{\alpha}(bt) &= \frac{1}{\pi bt\Gamma(\alpha)}\left(\frac{bt}{2}\right)^{\alpha} \int_0^{\frac{\pi}{2}} \int_0^{\pi}\biggr[2\sin^{2\alpha - 1}{\var}\sin{\phi}\biggr(\sin\left(bt\cos\var\sin{\phi}\right)
\\&\hspace{7cm} + \hspace{0.1cm} bt\cos\var\sin{\phi}\cos\left(bt\cos\var\sin{\phi}\right)\biggr]\biggr)\,\rmd{d}\phi \rmd{d}\var
\end{align*}

\newpage
\begin{align*}
\implies \int_0^{\infty} e^{-at}J_{\alpha}(bt)t^{-\alpha}\,\rmd{d}t &= \int_0^\infty \frac{1}{\pi b\Gamma(\alpha)}\left(\frac{b}{2}\right)^{\alpha} t^{-1}e^{-at}\Biggr(\int_0^{\frac{\pi}{2}} \int_0^{\pi}\left[2\sin^{2\alpha - 1}{\var}\sin{\phi}\left[\sin\left(bt\cos\var\sin{\phi}\right)\right.\right.
\\&\hspace{4cm}\left.\left.+ \hspace{0.1cm} bt\cos\var\sin{\phi}\cos\left(bt\cos\var\sin{\phi}\right)\right]\right]\,\rmd{d}\phi \rmd{d}\var\,\rmd{d}t\Biggr)
\\&= \frac{1}{\pi b\Gamma(\alpha)}\left(\frac{b}{2}\right)^{\alpha} \int_0^{\frac{\pi}{2}} \int_0^{\pi}\int_0^\infty t^{-1}e^{-at}\biggr(\left[2\sin^{2\alpha - 1}{\var}\sin{\phi}\left[\sin\left(bt\cos\var\sin{\phi}\right)\right.\right.
\\&\hspace{5.3cm}\left.\left.+ \hspace{0.1cm} bt\cos\var\sin{\phi}\cos\left(bt\cos\var\sin{\phi}\right)\right]\right]\,\rmd{d}t\rmd{d}\phi\,\rmd{d}\var\biggr)
\end{align*}\leavevmode\newline

From our previous calculation on  Lipschitz's result,
\begin{equation}\label{prev}
\int_0^{\pi}\left(\int_0^\infty e^{-at}\frac{\sin\left(bt\sin{\phi}\right)\sin{\phi}}{\pi bt}\rmd{d}t  + \int_0^\infty\frac{e^{-at}}{\pi}\sin^2{\phi}\cos\left(bt\sin{\phi}\right)\,\rmd{d}t\right)\rmd{d}\phi = \frac{1}{\sqrt{a^2 + b^2}}
\end{equation}

Replacing $b$ with $b\cos{\var}$ in \eqref{prev} above,
\begin{align*}
\implies &\int_0^{\pi}\int_0^\infty \frac{1}{b\pi\cos{\var}} t^{-1}e^{-at}\biggr(\sin{\phi}\left[\sin\left(bt\cos\var\sin{\phi}\right)+ bt\cos\var\sin{\phi}\cos\left(bt\cos\var\sin{\phi}\right)\right]\biggr)\,\rmd{d}t\rmd{d}\phi = \frac{1}{\sqrt{b^2\cos^2{\var} + a^2}} \\
\implies &\int_0^{\pi}\int_0^\infty \frac{1}{b\pi} t^{-1}e^{-at}\biggr(\sin{\phi}\left[\sin\left(bt\cos\var\sin{\phi}\right) + bt\cos\var\sin{\phi}\cos\left(bt\cos\var\sin{\phi}\right)\right]\biggr)\,\rmd{d}t\rmd{d}\phi = \frac{\cos{\var}}{\sqrt{b^2\cos^2{\var} + a^2}} 
\end{align*}
Hence,
\begin{align}
\int_0^{\infty} e^{-at}J_{\alpha}(bt)t^{-\alpha}\,\rmd{d}t  &= \frac{1}{\Gamma(\alpha)}\left(\frac{b}{2}\right)^\alpha \int_0^{\frac{\pi}{2}} \frac{2\sin^{2\alpha - 1}{\var}\cos{\var}}{\sqrt{b^2\cos^2{\var} + a^2}} \rmd{d}\var\nonumber
\\&= \frac{1}{\Gamma(\alpha)}\left(\frac{b}{2}\right)^\alpha\int_0^{\frac{\pi}{2}} \frac{2\sin^{2\alpha - 1}{\var}}{\sqrt{b^2 + a^2 - b^2\sin^2{\var}}} \rmd{d}\left(\sin{\var}\right)\nonumber
\\&= \frac{1}{\Gamma(\alpha)}\left(\frac{b}{2}\right)^\alpha\int_0^{1} \frac{2u^{2\alpha - 1}}{\sqrt{b^2 + a^2 - b^2u^2}} \rmd{d}u\nonumber
\\&= \frac{1}{\Gamma(\alpha)\sqrt{b^2 + a^2}}\left(\frac{b}{2}\right)^\alpha\int_0^{1} \frac{2u^{2\alpha - 1}}{\sqrt{1 - \frac{b^2u^2}{b^2 + a^2}}} \rmd{d}u \label{talkh}
\\&= \frac{1}{b\Gamma(\alpha)}\left(\frac{b}{2}\right)^\alpha\left(\frac{\sqrt{b^2 + a^2}}{b}\right)^{2\alpha - 1}\int_0^{\frac{b}{\sqrt{b^2 + a^2}}} \frac{2u^{2\alpha - 1}}{\sqrt{1 - u^2}} \rmd{d}u\nonumber
\\&= \frac{1}{2b\Gamma(\alpha)}\left(\frac{b}{2}\right)^\alpha\left(\frac{\sqrt{b^2 + a^2}}{b}\right)^{2\alpha - 1}\int_0^{\frac{b}{\sqrt{b^2 + a^2}}} \frac{2u^{\alpha - 1}}{\sqrt{1 - u}} \rmd{d}u\nonumber
\\&= \frac{1}{b\Gamma(\alpha)}\left(\frac{b}{2}\right)^\alpha\left(\frac{\sqrt{b^2 + a^2}}{b}\right)^{2\alpha - 1}\int_0^{\frac{b}{\sqrt{b^2 + a^2}}} u^{\alpha - 1} (1 - u)^{-\frac{1}{2}} \rmd{d}u\nonumber
\\&= \frac{1}{b\Gamma(\alpha)}\left(\frac{b}{2}\right)^\alpha\left(\frac{\sqrt{b^2 + a^2}}{b}\right)^{2\alpha - 1}B\left(\frac{b}{\sqrt{b^2 + a^2}}; \alpha, \frac{1}{2}\right)\nonumber
\\&= \frac{\left(b^2 + a^2\right)^\alpha}{2^{\alpha}b^{\alpha + 1}\Gamma(\alpha)}\cdot\frac{b}{\sqrt{b^2 + a^2}}B\left(\frac{b}{\sqrt{b^2 + a^2}}; \alpha, \frac{1}{2}\right)\label{ffo}\hspace{0.2cm}\,\,\, \textsf{valid for all}\,\,\, a \geq 0,\,\, \alpha > 0,
\end{align}
where $B(x; a, b)$ is the incomplete beta function.

\subsection{$\displaystyle\int_0^{\infty}J_{\frac{1}{2}}(bt)t^{-\frac{1}{2}}\,\rmd{d}t$}

From equation \eqref{talkh},
\begin{align*}
\int_0^{\infty} e^{-at}J_{\alpha}(bt)t^{-\alpha}\,\rmd{d}t &= \frac{1}{\Gamma(\alpha)\sqrt{b^2 + a^2}}\left(\frac{b}{2}\right)^\alpha\int_0^{1} \frac{2u^{2\alpha - 1}}{\sqrt{1 - \frac{b^2u^2}{b^2 + a^2}}} \rmd{d}u
\end{align*}

Let $\alpha = \frac{1}{2}$ and $a = 0$

\begin{align*}
\implies \int_0^{\infty} J_{\frac{1}{2}}(bt)t^{-\frac{1}{2}}\,\rmd{d}t &= \frac{1}{b\Gamma\left(\frac{1}{2}\right)}\left(\frac{b}{2}\right)^{\frac{1}{2}}\int_0^{1} \frac{2}{\sqrt{1 - u^2}} \rmd{d}u
\\&= \frac{2}{b\sqrt{\pi}}\sqrt{\frac{b}{2}} \arcsin{u}\vert_0^1
\\&= \frac{2}{b\sqrt{\pi}}\sqrt{\frac{b}{2}} \cdot \frac{\pi}{2}
\\&= \sqrt{\frac{\pi}{2b}}
\end{align*}

$$\therefore \int_0^{\infty}J_{\frac{1}{2}}(bt)t^{-\frac{1}{2}}\,\rmd{d}t = \sqrt{\frac{\pi}{2b}}\,\,\,\,\, \forall \,\,\Re(b) > 0$$

\subsection{$\displaystyle\int_0^{\infty}e^{-at} J_{\frac{1}{2}}(bt)t^{-\frac{1}{2}}\,\rmd{d}t$}

\begin{align*}
\int_0^{\infty}e^{-at} J_{\frac{1}{2}}(bt)t^{-\frac{1}{2}}\,\rmd{d}t &= \frac{1}{\Gamma\left(\frac{1}{2}\right)\sqrt{b^2 + a^2}}\left(\frac{b}{2}\right)^{\frac{1}{2}}\int_0^{1} \frac{2}{\sqrt{1 - \frac{b^2u^2}{b^2 + a^2}}} \rmd{d}u
\\&= \frac{2}{\sqrt{\pi}\sqrt{b^2 + a^2}}\left(\frac{b}{2}\right)^{\frac{1}{2}} \frac{\sqrt{b^2 + a^2}}{b}\arcsin\left(\frac{b}{\sqrt{b^2 + a^2}}u\right)\biggr\vert_0^1
\\&= \sqrt{\frac{2}{\pi b}} \arcsin\left(\frac{b}{\sqrt{b^2 + a^2}}\right)\,\,\,\,\, \forall \,\, \Re(b) > 0,\,\,\, \Re(a) \geq 0
\end{align*}

\subsection{$\displaystyle\int_0^{\infty}e^{-at} J_{1}(bt)t^{-1}\,\rmd{d}t$}

\begin{align*}
\int_0^{\infty} e^{-at}J_{1}(bt)t^{-1}\,\rmd{d}t &= \frac{1}{\Gamma(1)\sqrt{b^2 + a^2}}\left(\frac{b}{2}\right) \int_0^{1} \frac{2u}{\sqrt{1 - \frac{b^2u^2}{b^2 + a^2}}} \rmd{d}u
\\&= \frac{b}{2}\int_0^{1} \frac{2u}{\sqrt{b^2 + a^2 - b^2u^2}} \rmd{d}u
\\&=\frac{b}{2}\int_0^{1} \frac{1}{\sqrt{b^2 + a^2 - b^2u}} \rmd{d}u
\\&= \frac{b}{2}\cdot\frac{2}{-b^2}\sqrt{b^2 + a^2 - b^2u}\vert_0^1
\\&= -\frac{1}{b}\left(a - \sqrt{b^2 + a^2}\right) = \frac{1}{b}\left(\sqrt{b^2 + a^2} - a\right)\,\,\,\,\, \forall \,\, \Re(b) > 0,\,\,\, \Re(a) \geq 0
\end{align*}

\subsubsection{Special result}
$$\int_0^{\infty} e^{-3t}J_{1}(4t)t^{-1}\,\rmd{d}t = \frac{1}{4}\left(\sqrt{4^2 + 3^2} - 3\right) = \frac{1}{2}$$
$$\implies 2\int_0^{\infty} e^{-3t}J_{1}(4t)t^{-1}\,\rmd{d}t = 1$$

\subsection{$\displaystyle\int_0^{\infty}J_{1}(bt)t^{-1}\,\rmd{d}t$}

It has been established previously that
$$\int_0^{\infty} e^{-at}J_{1}(bt)t^{-1}\,\rmd{d}t  = \frac{1}{b}\left(\sqrt{b^2 + a^2} - a\right)$$

At $a=0$, 

$$\int_0^{\infty}J_{1}(bt)t^{-1}\,\rmd{d}t  = \frac{1}{b}\left(\sqrt{b^2 + 0^2} - 0\right) = 1$$

$$\therefore \int_0^{\infty}J_{1}(bt)t^{-1}\,\rmd{d}t = 1\,\,\,\,\, \forall \,\, \Re(b) > 0$$

\subsection{$\displaystyle\int_0^{\infty} e^{-at}J_{\alpha}(bt)t^{-\alpha}\,\rmd{d}t$}

From equation \eqref{ffo}
$$\int_0^{\infty} e^{-at}J_{\alpha}(bt)t^{-\alpha}\,\rmd{d}t =  \frac{\left(b^2 + a^2\right)^\alpha}{2^{\alpha}b^{\alpha + 1}\Gamma(\alpha)}\cdot\frac{b}{\sqrt{b^2 + a^2}}B\left(\frac{b}{\sqrt{b^2 + a^2}}; \alpha, \frac{1}{2}\right)\,\,\,\,\, \forall \,\, \alpha > 0,\,\,\, a \geq 0$$

At $a = 0$,
\begin{align*}
\int_0^{\infty}J_{\alpha}(bt)t^{-\alpha}\,\rmd{d}t &=  \frac{b^{\alpha - 1}}{2^{\alpha}\Gamma(\alpha)}B\left(1; \alpha, \frac{1}{2}\right)
\\&= \frac{b^{\alpha - 1}}{2^{\alpha}\Gamma(\alpha)}\beta\left(\alpha, \frac{1}{2}\right)
\\&= \frac{b^{\alpha - 1}}{2^{\alpha}\Gamma(\alpha)}\cdot\frac{\Gamma(\alpha)\Gamma\left(\frac{1}{2}\right)}{\Gamma\left(\alpha + \frac{1}{2}\right)}
\\&= \frac{b^{\alpha - 1}\sqrt{\pi}}{2^{\alpha}\Gamma\left(\alpha + \frac{1}{2}\right)}
\\&= \frac{b^{\alpha - 1}\sqrt{\pi}}{2^{\alpha}\left(\frac{\Gamma(2\alpha)\sqrt{\pi}}{\Gamma(\alpha)2^{2\alpha - 1}}\right)} = \frac{b^{\alpha - 1}}{\frac{\Gamma(2\alpha)}{\Gamma(\alpha)2^{\alpha - 1}}} 
\\&= \frac{\left(2b\right)^{\alpha - 1}\Gamma(\alpha)}{\Gamma(2\alpha)}
\end{align*}

$$\therefore \int_0^{\infty}J_{\alpha}(bt)t^{-\alpha}\,\rmd{d}t   =\begin{cases}
&\frac{\left(2b\right)^{\alpha - 1}\Gamma(\alpha)}{\Gamma(2\alpha)} \,\,\,\,\textsf{for $\alpha > 0$}
\\&\frac{b^{0 - 1}\sqrt{\pi}}{2^{0}\Gamma\left(0 + \frac{1}{2}\right)} = \frac{1}{b} \,\,\,\,\textsf{for $\alpha = 0$}
\end{cases}$$

\subsubsection{Special results}

$$\int_0^{\infty}J_{2}(3t)t^{-2}\,\rmd{d}t  = \frac{\left(2\times 3\right)^{1}\Gamma(2)}{\Gamma(4)} = 1$$
$$\int_0^{\infty}J_{0}(bt)\, \rmd{d}t  = \frac{1}{b}\,\,\,\,\,\forall \,\, \Re(b) > 0$$
$$\int_0^{\infty}J_{0}(t)\, \rmd{d}t  = 1\,\,\,\,\,\textsf{for $b = 1$}.$$

\section{Integral representations for $I_{\alpha}(x)$ and $I_{0}(x)$}

The Bessel functions are valid even for complex arguments $x$, and an important special case is that of a purely imaginary argument. In this case, the solutions to the Bessel equation are called the modified Bessel functions (or occasionally the hyperbolic Bessel functions) of the first and second kind and are defined as \cite{Treat5}

$$I_{\alpha}(x) = i^{-\alpha}J_{\alpha}(ix) = \sum_{m=0}^{\infty} \frac{1}{m!\Gamma(m + \alpha + 1)} \left(\frac{x}{2}\right)^{2m + \alpha}\hspace{0.5cm}\forall \,\, x \in \mathbb{C}$$

\begin{equation}\label{fin}
K_{\alpha}\left(x\right) = \frac{\pi}{2}\frac{I_{-\alpha}(x) - I_{\alpha}(x)}{\sin(\alpha \pi)}\hspace{1cm} \textsf{when $\alpha$ is not an integer.}
\end{equation}

Now, at $\alpha = 0$, 
\begin{equation}\label{what}
I_{0}(x) = J_{0}(ix)
\end{equation}

By equation \eqref{imag},
$$J_{0}(x) = \frac{1}{\pi x}\int_0^{\pi}\left(\sin\left(x\sin{\phi}\right) + x\sin{\phi}\cos\left(x\sin{\phi}\right)\right)\sin{\phi}\,\rmd{d}\phi$$

\begin{align*}
I_{0}(x) = J_{0}(ix) &= \frac{1}{\pi ix}\int_0^{\pi}\left(\sin\left(ix\sin{\phi}\right) + ix\sin{\phi}\cos\left(ix\sin{\phi}\right)\right)\sin{\phi}\,\rmd{d}\phi
\\&= \frac{1}{\pi x}\int_0^{\pi}\left(-i\sin\left(ix\sin{\phi}\right) + x\sin{\phi}\cosh\left(x\sin{\phi}\right)\right)\sin{\phi}\,\rmd{d}\phi
\\&= \frac{1}{\pi x}\int_0^{\pi}\left(\sinh\left(x\sin{\phi}\right) + x\sin{\phi}\cosh\left(x\sin{\phi}\right)\right)\sin{\phi}\,\rmd{d}\phi
\end{align*}

\begin{equation}\label{nam}
\implies I_{0}(x) = \frac{1}{\pi x}\int_0^{\pi}\left(\sinh\left(x\sin{\phi}\right) + x\sin{\phi}\cosh\left(x\sin{\phi}\right)\right)\sin{\phi}\,\rmd{d}\phi
\end{equation}\leavevmode\newline

By equation \eqref{refe},

$$J_{\alpha}(ix) = \frac{1}{\Gamma(\alpha)}\left(\frac{ix}{2}\right)^{\alpha} \int_0^1 J_{0}\left(ix\sqrt{1 - t}\right)t^{\alpha - 1} \rmd{d}t$$

$$\implies i^{-\alpha} J_{\alpha}(ix) = \frac{1}{\Gamma(\alpha)}\left(\frac{x}{2}\right)^{\alpha} \int_0^1 J_{0}\left(ix\sqrt{1 - t}\right)t^{\alpha - 1} \rmd{d}t$$\\

By equation \eqref{what}, $I_{0}(x) = J_{0}(ix)$.

$$\implies i^{-\alpha} J_{\alpha}(ix) = I_{\alpha}(x)= \frac{1}{\Gamma(\alpha)}\left(\frac{x}{2}\right)^{\alpha} \int_0^1 I_{0}\left(x\sqrt{1 - t}\right)t^{\alpha - 1} \rmd{d}t$$

\begin{equation}\label{name}
\implies  I_{\alpha}(x)= \frac{1}{\Gamma(\alpha)}\left(\frac{x}{2}\right)^{\alpha} \int_0^1 I_{0}\left(x\sqrt{1 - t}\right)t^{\alpha - 1} \rmd{d}t
\end{equation}

Inserting \eqref{nam} in \eqref{name},
\begin{align*}
I_{\alpha}(x) &= \frac{1}{\pi x\Gamma(\alpha)}\left(\frac{x}{2}\right)^{\alpha} \int_0^{\frac{\pi}{2}} \int_0^{\pi}\left[2\sin^{2\alpha - 1}{\var}\sin{\phi}\left[\sinh\left(x\cos\var\sin{\phi}\right) + x\cos\var\sin{\phi}\cosh\left(x\cos\var\sin{\phi}\right)\right]\right]\,\rmd{d}\phi \rmd{d}\var
\end{align*}

\section{Simplifications of the new integral representations}

\begin{align*}
J_{0}(x) &=\sum_{m=0}^{\infty} \frac{(-1)^m}{m!\Gamma(m + 1)} \left(\frac{x}{2}\right)^{2m}
\\&=\frac{1}{\pi x}\int_0^{\pi}\left(\sin\left(x\sin{\phi}\right) + x\sin{\phi}\cos\left(x\sin{\phi}\right)\right)\sin{\phi}\,\rmd{d}\phi
\\&= \frac{1}{\pi x}\int_0^{\pi}\left(\sin\left(x\sin{\phi}\right) + x\frac{\partial}{\partial x}\left(\sin\left(x\sin{\phi}\right)\right)\right)\sin{\phi}\,\rmd{d}\phi
\\&= \frac{1}{\pi x}\int_0^{\pi}\sin{\phi}\frac{\partial}{\partial x}\left(x\sin\left(x\sin{\phi}\right)\right)\rmd{d}\phi
\end{align*}

$$\implies J_0(x) = \frac{1}{\pi x}\int_0^{\pi}\sin{\phi}\frac{\partial}{\partial x}\left(x\sin\left(x\sin{\phi}\right)\right)\rmd{d}\phi$$

Since
$$I_0(x) = \frac{1}{\pi x}\int_0^{\pi}\left(\sinh\left(x\sin{\phi}\right) + x\sin{\phi}\cosh\left(x\sin{\phi}\right)\right)\sin{\phi}\,\rmd{d}\phi$$

$$\implies I_0(x) = \frac{1}{\pi x}\int_0^{\pi}\sin{\phi}\frac{\partial}{\partial x}\left(x\sinh\left(x\sin{\phi}\right)\right)\rmd{d}\phi$$

$$I_{\alpha}(x) = \frac{1}{\Gamma(\alpha)}\left(\frac{x}{2}\right)^{\alpha} \int_0^{\frac{\pi}{2}} 2\sin^{2\alpha - 1}{\var}\cos{\var}I_{0}\left(x\cos{\var}\right) \rmd{d}\var$$

\begin{align*}
\implies I_{\alpha}(x) &= \frac{1}{\pi x\Gamma(\alpha)}\left(\frac{x}{2}\right)^{\alpha} \int_0^{\frac{\pi}{2}} \int_0^\pi  2\sin^{2\alpha - 1}{\var}\cos{\var} \frac{\sin{\phi}}{\cos{\var}} \frac{\partial}{\partial \left(x\cos{\var}\right)}\left(x\cos{\var}\sinh\left(x\cos{\var}\sin{\phi}\right)\right) \rmd{d}\phi\rmd{d}\var
\\&= \frac{1}{\pi x\Gamma(\alpha)}\left(\frac{x}{2}\right)^{\alpha} \int_0^{\frac{\pi}{2}} \int_0^\pi 2\sin^{2\alpha - 1}{\var}\sin{\phi}\frac{\partial}{\partial \left(x\cos{\var}\right)}\left(x\cos{\var}\sinh\left(x\cos{\var}\sin{\phi}\right)\right) \rmd{d}\phi\rmd{d}\var
\end{align*}

Similarly,
$$J_{\alpha}(x) = \frac{1}{\pi x\Gamma(\alpha)}\left(\frac{x}{2}\right)^{\alpha} \int_0^{\frac{\pi}{2}} \int_0^\pi 2\sin^{2\alpha - 1}{\var}\sin{\phi}\frac{\partial}{\partial \left(x\cos{\var}\right)}\left(x\cos{\var}\sin\left(x\cos{\var}\sin{\phi}\right)\right) \rmd{d}\phi\rmd{d}\var$$

\section{An expression for $K_{\frac{1}{2}}\left(x\right)\,\,\, \forall \,\, \Re(x) > 0$ with the new integral representations}

$$I_{\alpha}(x) = \frac{1}{\Gamma(\alpha)}\left(\frac{x}{2}\right)^{\alpha} \int_0^{\frac{\pi}{2}} 2\sin^{2\alpha - 1}{\var}\cos{\var}I_{0}\left(x\cos{\var}\right) \rmd{d}\var$$

\begin{align*}
I_{\frac{1}{2}}(x) &= \frac{1}{\Gamma\left(\frac{1}{2}\right)}\left(\frac{x}{2}\right)^{\frac{1}{2}} \int_0^{\frac{\pi}{2}} 2\sin^{2\left(\frac{1}{2}\right) - 1}{\var}\cos{\var}I_{0}\left(x\cos{\var}\right) \rmd{d}\var
\\&= \frac{2}{\sqrt{\pi}}\left(\frac{x}{2}\right)^{\frac{1}{2}} \int_0^{\frac{\pi}{2}} \cos{\var}I_{0}\left(x\cos{\var}\right) \rmd{d}\var
\end{align*}

\begin{align*}
J_{-\frac{1}{2}}(x) &= \sum_{m=0}^{\infty} \frac{(-1)^m}{m!\Gamma(m - \frac{1}{2} + 1)} \left(\frac{x}{2}\right)^{2m - \frac{1}{2}}
\\&= \frac{1}{\Gamma\left(\frac{1}{2}\right)}\sqrt{\frac{2}{x}} +  \sum_{m=1}^{\infty} \frac{(-1)^m}{(m+1)!\Gamma(m - \frac{1}{2} + 1)} \left(\frac{x}{2}\right)^{2m - \frac{1}{2}}
\\&= \frac{1}{\Gamma\left(\frac{1}{2}\right)}\sqrt{\frac{2}{x}} +  \sum_{m=0}^{\infty} \frac{(-1)^{m+1}}{(m+1)!\Gamma(m + \frac{1}{2} + 1)} \left(\frac{x}{2}\right)^{2m + \frac{3}{2}}
\\&= \frac{1}{\Gamma\left(\frac{1}{2}\right)}\sqrt{\frac{2}{x}} - \sum_{m=0}^{\infty} \frac{(-1)^m\Gamma(m + 1)\Gamma\left(\frac{1}{2}\right)}{\Gamma(m + 2)\Gamma(m + 1)\Gamma(m + \frac{1}{2} + 1)} \left(\frac{x}{2}\right)^{2m + \frac{3}{2}} 
\\&= \frac{1}{\Gamma\left(\frac{1}{2}\right)}\sqrt{\frac{2}{x}} - \frac{1}{\Gamma\left(\frac{1}{2}\right)}\sum_{m=0}^{\infty} \frac{(-1)^m\Gamma(m + 1)\Gamma\left(\frac{1}{2}\right)}{m!\Gamma(m + 2)\Gamma(m + \frac{1}{2} + 1)} \left(\frac{x}{2}\right)^{2m + \frac{3}{2}}
\\&= \frac{1}{\Gamma\left(\frac{1}{2}\right)}\sqrt{\frac{2}{x}} - \frac{1}{\Gamma\left(\frac{1}{2}\right)}\sum_{m=0}^{\infty} \frac{(-1)^m}{m!\Gamma(m + 2)}\beta\left(m + 1, \frac{1}{2}\right) \left(\frac{x}{2}\right)^{2m + \frac{3}{2}}
\\&= \frac{1}{\Gamma\left(\frac{1}{2}\right)}\sqrt{\frac{2}{x}} -  \frac{1}{\Gamma\left(\frac{1}{2}\right)}\sum_{m=0}^{\infty} \frac{(-1)^m}{m!\Gamma(m + 2)} \int_0^1 t^{m} \left(1 - t\right)^{\frac{1}{2} - 1} \rmd{d}t \left(\frac{x}{2}\right)^{2m + \frac{3}{2}}
\\&= \frac{1}{\Gamma\left(\frac{1}{2}\right)}\sqrt{\frac{2}{x}} - \frac{1}{\Gamma\left(\frac{1}{2}\right)}\left(\frac{x}{2}\right)^{\frac{1}{2}} \int_0^1 \sum_{m=0}^{\infty} \frac{(-1)^m}{m!\Gamma(m + 2)} t^{m}  \left(\frac{x}{2}\right)^{2m + 1} \left(1 - t\right)^{-\frac{1}{2}} \rmd{d}t
\\&= \frac{1}{\Gamma\left(\frac{1}{2}\right)}\sqrt{\frac{2}{x}} - \frac{1}{\Gamma\left(\frac{1}{2}\right)}\left(\frac{x}{2}\right)^{\frac{1}{2}} \int_0^1 \sum_{m=0}^{\infty} \frac{(-1)^m}{m!\Gamma(m + 2)}  \left(\frac{x\sqrt{t}}{2}\right)^{2m + 1} t^{-\frac{1}{2}}\left(1 - t\right)^{-\frac{1}{2}} \rmd{d}t
\end{align*}

\begin{align*}
\implies J_{-\frac{1}{2}}(x) &=  \frac{1}{\Gamma\left(\frac{1}{2}\right)}\sqrt{\frac{2}{x}} - \frac{1}{\Gamma\left(\frac{1}{2}\right)}\left(\frac{x}{2}\right)^{\frac{1}{2}} \int_0^1 J_{1}(x\sqrt{t})t^{-\frac{1}{2}}\left(1 - t\right)^{-\frac{1}{2}} \rmd{d}t
\\&=  \frac{1}{\Gamma\left(\frac{1}{2}\right)}\sqrt{\frac{2}{x}} - \frac{1}{\Gamma\left(\frac{1}{2}\right)}\left(\frac{x}{2}\right)^{\frac{1}{2}} \int_0^1 J_{1}\left(x\sqrt{1 - t}\right)t^{\frac{1}{2} - 1}\left(1 - t\right)^{-\frac{1}{2}} \rmd{d}t
\end{align*}\leavevmode\newline

Since, $I_{\alpha}(x) = i^{-\alpha}J_{\alpha}(ix)$, 
$$\implies I_{-\frac{1}{2}}(x) = i^{\frac{1}{2}}J_{\frac{1}{2}}(ix)$$

$$I_{1}(x) = i^{-1}J_{1}(ix)\implies J_{1}(ix) = iI_{1}(x) \implies iJ_{1}(ix) = -I_{1}(x)$$

\begin{align*}
\implies i^{\frac{1}{2}}J_{\frac{1}{2}}(ix) &=  \frac{1}{\Gamma\left(\frac{1}{2}\right)}\sqrt{\frac{2}{x}} - \frac{i}{\Gamma\left(\frac{1}{2}\right)}\left(\frac{x}{2}\right)^{\frac{1}{2}} \int_0^1 J_{1}\left(ix\sqrt{1 - t}\right)t^{\frac{1}{2} - 1}\left(1 - t\right)^{-\frac{1}{2}} \rmd{d}t
\\& \implies I_{-\frac{1}{2}}(x) = \frac{1}{\Gamma\left(\frac{1}{2}\right)}\sqrt{\frac{2}{x}} + \frac{1}{\Gamma\left(\frac{1}{2}\right)}\left(\frac{x}{2}\right)^{\frac{1}{2}} \int_0^1 I_{1}\left(x\sqrt{1 - t}\right)t^{-\frac{1}{2}}\left(1 - t\right)^{-\frac{1}{2}} \rmd{d}t
\end{align*}\leavevmode\newline

\begin{align*}
I_{-\frac{1}{2}}(x) &=  \frac{1}{\Gamma\left(\frac{1}{2}\right)}\sqrt{\frac{2}{x}} +  \frac{1}{\Gamma\left(\frac{1}{2}\right)}\left(\frac{x}{2}\right)^{\frac{1}{2}} \int_0^{\frac{\pi}{2}} 2\sin^{2\left(\frac{1}{2}\right) - 1}{\var}\cos{\var}\frac{1}{\cos{\var}}I_{1}\left(x\cos{\var}\right) \rmd{d}\var
\\&= \frac{1}{\Gamma\left(\frac{1}{2}\right)}\sqrt{\frac{2}{x}} + \frac{2}{\Gamma\left(\frac{1}{2}\right)}\left(\frac{x}{2}\right)^{\frac{1}{2}} \int_0^{\frac{\pi}{2}} I_{1}\left(x\cos{\var}\right) \rmd{d}\var
\end{align*}

\begin{align*}
\int_0^{\frac{\pi}{2}} \cos{\var} I_0\left(x\cos{\var}\right) \rmd{d}\var &= \int_0^1 I_0\left(x\sqrt{1 - u^2}\right)\rmd{d}u
\\&= \frac{1}{2}\int_0^1 u^{-\frac{1}{2}} I_0\left(x\sqrt{1 - u}\right)\rmd{d}u
\\&= \frac{1}{2}\int_0^1 \left(1 - u\right)^{-\frac{1}{2}} I_0\left(x\sqrt{u}\right)\rmd{d}u
\\&= \frac{1}{2}\sum_{m=0}^{\infty} \frac{1}{m!\Gamma(m + 1)} \int_0^1 \left(1 - u\right)^{-\frac{1}{2}} \left(\frac{x\sqrt{u}}{2}\right)^{2m} \,\, \rmd{d}u
\\&= \frac{1}{2}\sum_{m=0}^{\infty} \frac{1}{m!\Gamma(m + 1)}\left(\frac{x}{2}\right)^{2m} \int_0^1 u^m \left(1 - u\right)^{-\frac{1}{2}}  \,\, \rmd{d}u
\\&= \frac{1}{2}\sum_{m=0}^{\infty} \frac{1}{m!\Gamma(m + 1)}\left(\frac{x}{2}\right)^{2m}\beta\left(m + 1, \frac{1}{2}\right)
\\&= \frac{1}{2}\sum_{m=0}^{\infty} \frac{1}{m!\Gamma(m + 1)}\left(\frac{x}{2}\right)^{2m}\frac{\Gamma(m + 1)\Gamma\left(\frac{1}{2}\right)}{\Gamma\left(m + \frac{3}{2}\right)}
\\&= \frac{1}{2}\sum_{m=0}^{\infty} \frac{1}{m!}\left(\frac{x}{2}\right)^{2m}\frac{\Gamma\left(\frac{1}{2}\right)}{\Gamma\left(m + \frac{3}{2}\right)}
\\&= \frac{1}{2}\sum_{m=0}^{\infty} \frac{1}{\Gamma(m + 1)}\left(\frac{x}{2}\right)^{2m}\frac{\Gamma\left(\frac{1}{2}\right)}{\Gamma\left(m + \frac{3}{2}\right)}
\end{align*}\leavevmode\newline

By Legendre's duplication formula,
$$\frac{2^{2m + 1}}{\Gamma(2m + 2)} = \frac{\Gamma\left(\frac{1}{2}\right)}{\Gamma\left(m + \frac{3}{2}\right)\Gamma(m + 1)}$$

\begin{align*}
\implies \int_0^{\frac{\pi}{2}} \cos{\var} I_0\left(x\cos{\var}\right) \rmd{d}\var &= \frac{1}{2}\sum_{m=0}^{\infty} \frac{2^{2m + 1}}{\Gamma(2m + 2)} \left(\frac{x}{2}\right)^{2m}
\\&= \sum_{m=0}^{\infty} \frac{x^{2m}}{\Gamma(2m + 2)} = \sum_{m=0}^{\infty} \frac{x^{2m}}{(2m + 1)!}
\\&= \frac{1}{x} \sum_{m=0}^{\infty} \frac{x^{2m + 1}}{(2m + 1)!}
\\&= \frac{\sinh(x)}{x}
\end{align*}

\begin{align*}
\implies I_{\frac{1}{2}}(x) &= \frac{2}{\sqrt{\pi}}\left(\frac{x}{2}\right)^{\frac{1}{2}} \int_0^{\frac{\pi}{2}} \cos{\var}I_{0}\left(x\cos{\var}\right) \rmd{d}\var
\\&= \sqrt{\frac{2x}{\pi}}\frac{\sinh(x)}{x}
\\&=\sqrt{\frac{2}{\pi x}}\sinh(x)
\end{align*}

\begin{align*}
\int_0^{\frac{\pi}{2}} I_1\left(x\cos{\var}\right) \rmd{d}\var &= \int_0^1 \frac{I_1\left(x\sqrt{1 - u^2}\right)}{\sqrt{1 - u^2}}\rmd{d}u
\\&= \frac{1}{2}\int_0^1 u^{-\frac{1}{2}} \left(1 - u\right)^{-\frac{1}{2}} I_1\left(x\sqrt{1 - u}\right)\rmd{d}u
\\&= \frac{1}{2}\int_0^1 u^{-\frac{1}{2}}  \left(1 - u\right)^{-\frac{1}{2}} I_1\left(x\sqrt{u}\right)\rmd{d}u
\\&= \frac{1}{2}\sum_{m=0}^{\infty} \frac{1}{m!\Gamma(m + 2)} \int_0^1 u^{-\frac{1}{2}} \left(1 - u\right)^{-\frac{1}{2}} \left(\frac{x\sqrt{u}}{2}\right)^{2m + 1} \,\, \rmd{d}u
\\&= \frac{1}{2}\sum_{m=0}^{\infty} \frac{1}{m!\Gamma(m + 2)}\left(\frac{x}{2}\right)^{2m + 1} \int_0^1 u^{m} \left(1 - u\right)^{-\frac{1}{2}}  \,\, \rmd{d}u
\\&= \frac{1}{2}\sum_{m=0}^{\infty} \frac{1}{m!\Gamma(m + 2)}\left(\frac{x}{2}\right)^{2m + 1}\beta\left(m+1, \frac{1}{2}\right)
\\&= \frac{1}{2}\sum_{m=0}^{\infty} \frac{1}{m!(m + 1)\Gamma(m + 1)}\left(\frac{x}{2}\right)^{2m + 1}\frac{\Gamma(m + 1)\Gamma\left(\frac{1}{2}\right)}{\Gamma\left(m + \frac{3}{2}\right)}
\\&= \frac{1}{2}\sum_{m=0}^{\infty} \frac{1}{m!(m + 1)}\left(\frac{x}{2}\right)^{2m + 1}\frac{\Gamma\left(\frac{1}{2}\right)}{\Gamma\left(m + \frac{3}{2}\right)}
\\&= \frac{1}{2}\sum_{m=0}^{\infty} \frac{1}{\Gamma(m + 2)}\left(\frac{x}{2}\right)^{2m + 1}\frac{\Gamma\left(\frac{1}{2}\right)}{\Gamma\left(m + \frac{3}{2}\right)}
\\&= \frac{1}{2}\sum_{m=0}^{\infty} \frac{m + \frac{3}{2}}{\Gamma(m + 2)}\left(\frac{x}{2}\right)^{2m + 1}\frac{\Gamma\left(\frac{1}{2}\right)}{\Gamma\left(m + \frac{5}{2}\right)}
\\&= \frac{1}{4}\sum_{m=0}^{\infty} \frac{2m + 3}{\Gamma(m + 2)}\left(\frac{x}{2}\right)^{2m + 1}\frac{\Gamma\left(\frac{1}{2}\right)}{\Gamma\left(m + \frac{5}{2}\right)}
\end{align*}

\begin{align*}
\implies \int_0^{\frac{\pi}{2}} \cos{\var} I_0\left(x\cos{\var}\right) \rmd{d}\var &= \frac{1}{4}\sum_{m=0}^{\infty} \frac{2^{2m + 3}}{\Gamma(2m + 4)}\left(2m + 3\right) \left(\frac{x}{2}\right)^{2m + 1}
\\&= 2\cdot\frac{x}{2}\sum_{m=0}^{\infty} \frac{x^{2m}\left(2m + 3\right)}{(2m + 3)\Gamma(2m + 3)} = x\sum_{m=0}^{\infty} \frac{x^{2m}}{(2m + 2)!}
\\&= \frac{1}{x} \sum_{m=0}^{\infty} \frac{x^{2m + 2}}{(2m + 2)!}
\\&= \frac{\cosh(x) - 1}{x}
\end{align*}

\begin{align*}
I_{-\frac{1}{2}}(x) &= \frac{1}{\Gamma\left(\frac{1}{2}\right)}\sqrt{\frac{2}{x}} + \sqrt{\frac{2x}{\pi}}\left(\frac{\cosh(x) - 1}{x}\right)
\\&= \sqrt{\frac{2}{\pi x}}+ \sqrt{\frac{2x}{\pi}}\left(\frac{\cosh(x) - 1}{x}\right)
\\&= \sqrt{\frac{2}{\pi x}}\cosh(x)
\end{align*}\leavevmode\newline

\begin{align*}
I_{-\frac{1}{2}}(x) -  I_{\frac{1}{2}}(x) &= \sqrt{\frac{2}{\pi x}}\cosh(x) - \sqrt{\frac{2}{\pi x}}\sinh(x)
\\&= \sqrt{\frac{2}{\pi x}} e^{-x}
\end{align*}

At $\alpha = \frac{1}{2}$ in equation \eqref{fin},
\begin{align*}
K_{\frac{1}{2}}\left(x\right) &= \frac{\pi}{2}\frac{I_{-\frac{1}{2}}(x) - I_{\frac{1}{2}}(x)}{\sin\left(\frac{\pi}{2}\right)}
\\&= \frac{\pi}{2}\left(I_{-\frac{1}{2}}(x) - I_{\frac{1}{2}}\right) 
\\&= \frac{\pi}{2} \cdot \sqrt{\frac{2}{\pi x}} e^{-x} = \sqrt{\frac{\pi}{2 x}} e^{-x}
\end{align*}

\begin{equation}\label{speci}
\therefore K_{\frac{1}{2}}\left(x\right)  = \sqrt{\frac{\pi}{2 x}} e^{-x}
\end{equation}

\section{Deriving another integral representation for $K_{\alpha}(xz)$ directly from Basset's formula}

By Basset's formula \cite{Treat4},
$$K_{\alpha}(xz) = \frac{\Gamma\left(\alpha + \frac{1}{2}\right)\cdot (2z)^\alpha}{x^\alpha \Gamma\left(\frac{1}{2}\right)} \int_0^{\infty} \frac{\cos(xu)}{\left(u^2 + z^2\right)^{\alpha + \frac{1}{2}}} \rmd{d}u$$

By elmentary substitution, Basset's formula can be rewritten to $K_{\alpha}(z).$
\begin{align*}
K_{\alpha}(xz) &= \frac{\Gamma\left(\alpha + \frac{1}{2}\right)\cdot (2z)^\alpha}{x^\alpha \Gamma\left(\frac{1}{2}\right)} \int_0^{\infty} \frac{\cos(xu)}{\left(u^2 + z^2\right)^{\alpha + \frac{1}{2}}} \rmd{d}u
\\&= \frac{\Gamma\left(\alpha + \frac{1}{2}\right)\cdot (2z)^\alpha}{x^\alpha \Gamma\left(\frac{1}{2}\right)} z\int_0^{\infty} \frac{\cos(xzu)}{\left(z^2 u^2 + z^2\right)^{\alpha + \frac{1}{2}}} \rmd{d}u
\\&= \frac{\Gamma\left(\alpha + \frac{1}{2}\right)\cdot (2z)^\alpha}{(z^{2\alpha}\cdot z)x^\alpha \Gamma\left(\frac{1}{2}\right)} z\int_0^{\infty} \frac{\cos(xzu)}{\left(u^2 + 1\right)^{\alpha + \frac{1}{2}}} \rmd{d}u
\\&= \frac{\Gamma\left(\alpha + \frac{1}{2}\right)\cdot (2)^\alpha}{z^{\alpha}x^\alpha \Gamma\left(\frac{1}{2}\right)} \int_0^{\infty} \frac{\cos(xzu)}{\left(u^2 + 1\right)^{\alpha + \frac{1}{2}}} \rmd{d}u
\\\implies K_{\alpha}(z)&= \frac{\Gamma\left(\alpha + \frac{1}{2}\right)\cdot (2)^\alpha}{z^{\alpha}\Gamma\left(\frac{1}{2}\right)} \int_0^{\infty} \frac{\cos(zu)}{\left(u^2 + 1\right)^{\alpha + \frac{1}{2}}} \rmd{d}u
\end{align*}

Consider the following integral
$$H(\beta, p, R) =  \frac{2\Gamma(p)}{\sqrt{\pi}} \int_0^\infty \frac{\cos\left(2Rx\right)}{\left(\beta^2 + x^2\right)^p} \rmd{d}x$$

\begin{align*}
H(\beta, p, R) &= \frac{2\Gamma(p)}{\sqrt{\pi}} \int_0^\infty \frac{\cos\left(2Rx\right)}{\left(\beta^2 + x^2\right)^p} \rmd{d}x
\\&=  \frac{2}{\sqrt{\pi}} \int_0^\infty\int_0^\infty  t^{p - 1}e^{-\left(\beta^2 + x^2\right)t} \cos\left(2Rx\right) \rmd{d}t\,\rmd{d}x
\\&=  \frac{2}{\sqrt{\pi}} \int_0^\infty\int_0^\infty  t^{p - 1}e^{-\left(\beta^2 + x^2\right)t} \cos\left(2Rx\right) \rmd{d}x\,\rmd{d}t
\\&=  \frac{2}{\sqrt{\pi}} \int_0^\infty t^{p - 1}e^{-\beta^2 t} \int_0^\infty e^{-x^2t} \cos\left(2Rx\right) \rmd{d}x\,\rmd{d}t
\end{align*}

Let
$$f(R, t) = \int_0^\infty e^{-x^2t} \cos\left(2Rx\right) \rmd{d}x$$

\begin{align*}
\frac{\partial f(R, t)}{\partial R} &= -2\int_0^\infty xe^{-x^2t} \sin\left(2Rx\right) \rmd{d}x
\\&= \frac{1}{t}\int_0^\infty \sin\left(2Rx\right) \rmd{d}\left(e^{-x^2t}\right)
\\&= -\frac{2R}{t}\int_0^\infty e^{-x^2t} \cos\left(2Rx\right) \rmd{d}x = -\frac{2R}{t}f(R, t) 
\end{align*}

$$\frac{\partial f(R, t)}{\partial R} = -\frac{2R}{t}f(R, t)$$
$$\implies \frac{1}{f(R, t)}\frac{\partial f(R, t)}{\partial R} = -\frac{2R}{t}$$
$$\implies \frac{\partial \left(\ln\left(f(R, t)\right)\right)}{\partial R} = -\frac{2R}{t}$$
$$\implies \partial \left(\ln\left(f(R, t)\right)\right)= -\frac{2R}{t}\partial R$$
$$\implies \int \rmd{d}\left(\ln\left(f(R, t)\right)\right)= -\int \frac{2R}{t}\rmd{d}R$$
$$\implies \ln\left(f(R, t)\right) = -\frac{R^2}{t} + C(t)$$
$$\implies f(R, t) = A(t)e^{-\frac{R^2}{t}}$$

$$\implies f(R, t) = \int_0^\infty e^{-x^2t} \cos\left(2Rx\right) \rmd{d}x = A(t)e^{-\frac{R^2}{t}}$$

$$\therefore f(0, t) = \int_0^\infty e^{-x^2t} \rmd{d}x = A(t).$$

\begin{align*}
\int_0^\infty e^{-x^2t} \rmd{d}x &\stackrel{x \rightarrow \sqrt{u}}{=} \frac{1}{2} \int_0^{\infty} u^{-\frac{1}{2}} e^{-ut}\,\,\rmd{d}u
\\&= \frac{1}{2}\mathcal{L}\left(u^{-\frac{1}{2}}\right)(t) = \frac{\Gamma\left(\frac{1}{2}\right)}{2t^{\frac{1}{2}}} = \frac{1}{2}\sqrt{\frac{\pi}{t}}
\end{align*}

$$\implies f(R, t) = \int_0^\infty e^{-x^2t} \cos\left(2Rx\right) \rmd{d}x = \frac{1}{2}\sqrt{\frac{\pi}{t}} e^{-\frac{R^2}{t}}$$

\begin{align*}
\implies H(\beta, p, R) &= \frac{2}{\sqrt{\pi}} \int_0^\infty t^{p - 1}e^{-\beta^2 t} \int_0^\infty e^{-x^2t} \cos\left(2Rx\right) \rmd{d}x\,\rmd{d}t
\\&= \frac{2}{\sqrt{\pi}} \int_0^\infty t^{p - 1}e^{-\beta^2 t}\left(\frac{1}{2}\sqrt{\frac{\pi}{t}} e^{-\frac{R^2}{t}}\right)\rmd{d}t
\\&=\int_0^{\infty} t^{p - \frac{3}{2}} e^{-\left(\beta^2 t + \frac{R^2}{t}\right)} \rmd{d}t
\end{align*}

$$\implies H(\beta, p, R) = \int_0^{\infty} t^{p - \frac{3}{2}} e^{-\left(\beta^2 t + \frac{R^2}{t}\right)} \rmd{d}t = \frac{2\Gamma(p)}{\sqrt{\pi}} \int_0^\infty \frac{\cos\left(2Rt\right)}{\left(\beta^2 + t^2\right)^p} \rmd{d}t$$

Replacing $\beta^2$  with $a$, $R^2$ with $z$,
$$\int_0^{\infty} t^{p - \frac{3}{2}} e^{-\left(a t + \frac{z}{t}\right)} \rmd{d}t = \frac{2\Gamma(p)}{\sqrt{\pi}} \int_0^\infty \frac{\cos\left(2\sqrt{z}t\right)}{\left(a^2 + t^2\right)^p} \rmd{d}t$$

Replacing $p$ with $\frac{1}{2} + \alpha$, $z$ with $\frac{z^2}{4}$, $a$ with 1,
$$\int_0^{\infty} t^{\alpha - 1} e^{-\left( t + \frac{z^2}{4t}\right)} \rmd{d}t = \frac{2\Gamma\left(\frac{1}{2} + \alpha\right)}{\sqrt{\pi}} \int_0^\infty \frac{\cos\left(zt\right)}{\left(1 + t^2\right)^{\frac{1}{2} + \alpha}} \rmd{d}t$$

Multiplying both sides by $\frac{2^{\alpha - 1}}{z^{\alpha}}$,
$$\frac{\Gamma\left(\alpha + \frac{1}{2}\right)\cdot (2)^\alpha}{z^{\alpha}\Gamma\left(\frac{1}{2}\right)} \int_0^{\infty} \frac{\cos(zt)}{\left(t^2 + 1\right)^{\alpha + \frac{1}{2}}} \rmd{d}t = \frac{1}{z^{\alpha}} \int_0^{\infty} \left(2t\right)^{\alpha - 1} e^{-\left(t + \frac{z^2}{4t}\right)} \rmd{d}t$$

Thus if 
$$K_{\alpha}(z) = \frac{\Gamma\left(\alpha + \frac{1}{2}\right)\cdot (2)^\alpha}{z^{\alpha}\Gamma\left(\frac{1}{2}\right)} \int_0^{\infty} \frac{\cos(zt)}{\left(t^2 + 1\right)^{\alpha + \frac{1}{2}}} \rmd{d}t$$

Then
$$K_{\alpha}(z) = \frac{1}{z^{\alpha}} \int_0^{\infty} \left(2t\right)^{\alpha - 1} e^{-\left(t + \frac{z^2}{4t}\right)} \rmd{d}t,$$

$$\implies K_{\alpha}(z) = \frac{\Gamma\left(\alpha + \frac{1}{2}\right)\cdot (2)^\alpha}{z^{\alpha}\Gamma\left(\frac{1}{2}\right)} \int_0^{\infty} \frac{\cos(zt)}{\left(t^2 + 1\right)^{\alpha + \frac{1}{2}}} \rmd{d}t = \frac{1}{z^{\alpha}} \int_0^{\infty} \left(2t\right)^{\alpha - 1} e^{-\left(t + \frac{z^2}{4t}\right)} \rmd{d}t$$

\section{Another variant of Hardy's integral in relation to partial differential equations}

Prof. Godfrey Harold Hardy, FRS gave the following integral reprsentation \cite{Treat1}
$$\int_0^\infty \sin\left(au + \frac{b}{u}\right) \frac{\rmd{d}u}{u} = \pi J_{0}\left\{2\sqrt{ab}\right\}$$

On the transformation $u \longrightarrow u^2$, this is equivalent to
$$\int_0^\infty \sin\left(au^2 + \frac{b}{u^2}\right) \frac{\rmd{d}u}{u} = \frac{\pi}{2} J_{0}\left\{2\sqrt{ab}\right\}$$

Sequel to the above, a similar integral representation for $K_{\frac{1}{2}}\left(z\right)$ can be proposed in what follows
\begin{proposition*}
$$\int_0^\infty \sin\left(au^2 - \frac{b}{u^2}\right) \rmd{d}u = \left(\frac{b}{2a}\right)^{\frac{1}{4}} K_{\frac{1}{2}}\left\{2\sqrt{ab}\right\}$$
\end{proposition*}

\begin{proof}
We can take the partial derivatives of the integrand wrt. $a$ or $b$ by differentiating underthe integral sign since the partial derivatives are continuous in the limits of the integral.

$$\frac{\partial}{\partial a}\left(\int_0^\infty \sin\left(au^2 - \frac{b}{u^2}\right) \rmd{d}u \right) = \int_0^\infty u^2 \cos\left(au^2 - \frac{b}{u^2}\right) \rmd{d}u$$

\begin{align*}
\frac{\partial^2}{\partial a\partial b}\left(\int_0^\infty \sin\left(au^2 - \frac{b}{u^2}\right) \rmd{d}u \right) &= \int_0^\infty -\frac{1}{u^2} \cdot u^2 -\sin\left(au^2 - \frac{b}{u^2}\right) \rmd{d}u
\\&= \int_0^\infty \sin\left(au^2 - \frac{b}{u^2}\right) \rmd{d}u 
\end{align*}

$$\frac{\partial^2}{\partial a\partial b}\left(\int_0^\infty \sin\left(au^2 - \frac{b}{u^2}\right) \rmd{d}u \right) = \int_0^\infty \sin\left(au^2 - \frac{b}{u^2}\right) \rmd{d}u $$

Let 
$$\int_0^\infty \sin\left(au^2 - \frac{b}{u^2}\right) \rmd{d}u = F(a, b)$$

$$\implies \frac{\partial^2}{\partial a\partial b}F(a, b) = F(a, b)$$

Let $F(a, b)$ be the product of two functions $A(a)$, $B(b)$,
$$\implies F(a, b) = A(a)B(b)$$
$$\implies A'(a)B'(b) = A(a)B(b)$$
$$\implies \frac{A'(a)}{A(a)} = \frac{B(b)}{B'(b)}$$

Let 
$$\frac{A'(a)}{A(a)} = \frac{B(b)}{B'(b)} = \lambda(a, b)$$

$$\implies \frac{A'(a)}{A(a)} = \frac{B(b)}{B'(b)} = \lambda(a, b)$$

$$\implies \frac{A'(a)}{A(a)} = \lambda(a, b),\hspace{0.2cm} \frac{B(b)}{B'(b)} = \lambda(a, b)$$
$$\implies \frac{A'(a)}{A(a)} = \lambda(a, b),\hspace{0.2cm} \frac{B'(b)}{B(b)} = \frac{1}{\lambda(a, b)}$$

$$\implies \frac{\rmd{d}\left(A(a)\right)}{A(a)} = \lambda(a, b)\,\,\rmd{d}a, \hspace{0.2cm} \frac{\rmd{d}\left(B(b)\right)}{B(b)} = \frac{1}{\lambda(a, b)}\rmd{d}b$$
$$\implies \int \frac{\rmd{d}\left(A(a)\right)}{A(a)} = \int \lambda(a, b)\,\,\rmd{d}a,  \hspace{0.2cm} \int \frac{\rmd{d}\left(B(b)\right)}{B(b)} = \int \frac{1}{\lambda(a, b)}\rmd{d}b$$
$$\implies \ln\left(A(a)\right) = \int \lambda(a, b)\,\,\rmd{d}a + c_1, \hspace{0.2cm} \ln\left(B(b)\right) =\int \frac{\rmd{d}b}{\lambda(a, b)} + c_2$$
$$\implies A(a) = C_1 e^{\int \lambda(a, b)\,\,\rmd{d}a} , \hspace{0.2cm} B(b) = C_2 e^{\int \frac{\rmd{d}b}{\lambda(a, b)}}$$\\

Since $F(a, b) = A(a)B(b)$
$$\implies F(a, b) = C_1 C_2 e^{\int \lambda(a, b)\,\,\rmd{d}a + \int \frac{\rmd{d}b}{\lambda(a, b)}} = C_3 \exp\left(\int \lambda(a, b)\,\,\rmd{d}a + \int \frac{\rmd{d}b}{\lambda(a, b)}\right)$$

\begin{align*}
F(a, 0) &= C_3 \exp\left(\int \lambda(a, 0)\,\,\rmd{d}a + \int \frac{\rmd{d}b}{\lambda(a, 0)}\right)
\\&= \int_0^\infty \sin\left(au^2 - \frac{0}{u^2}\right) \rmd{d}u = \int_0^\infty \sin\left(au^2\right) \rmd{d}u
\\&= \frac{1}{2\sqrt{a}}\sqrt{\frac{\pi}{2}} \exp(0)
\end{align*}

$$\implies C_3 = \frac{1}{2\sqrt{a}}\sqrt{\frac{\pi}{2}},\hspace{0.5cm} \int \lambda(a, 0)\,\,\rmd{d}a + \int \frac{\rmd{d}b}{\lambda(a, 0)} = 0$$

Let $\lambda(a, 0) = f(a)$

$$\int f(a)\,\,\rmd{d}a + \int \frac{\rmd{d}b}{f(a)} = 0$$
$$\int f(a)\,\,\rmd{d}a = -\int \frac{\rmd{d}b}{f(a)} = -\frac{b}{f(a)} + \textup{constant}$$

Taking the derivative at $a$,
$$\frac{\rmd{d}}{\rmd{d}a}\left(\int f(a)\,\,\rmd{d}a\right) = -\frac{\rmd{d}}{\rmd{d}a}\left(\frac{b}{f(a)}\right)$$

$$\implies f(a) = -b\cdot\frac{-f'(a)}{f(a)^2} = \frac{bf'(a)}{f(a)^2}$$

$$\implies f'(a) = \frac{f(a)^3}{b}$$
$$\implies \frac{\rmd{d}\left(f(a)\right)}{f(a)^3} = \frac{\rmd{d}a}{b}$$
$$\implies \int \frac{\rmd{d}\left(f(a)\right)}{f(a)^3} = \int \frac{\rmd{d}a}{b}$$

For $\lambda_1, \gamma_1$ constant,
$$\implies -\frac{1}{2f(a)^2} = \frac{a}{b} + \lambda_1$$
$$\implies \frac{1}{f(a)^2} = -\frac{2a}{b} + \lambda_1 = \frac{b\gamma_1 - 2a}{b}$$
$$\implies f(a) = \sqrt{\frac{b}{b\gamma_1 - 2a}}$$

\begin{align*}
F(0, b) &= C_3 \exp\left(\int \lambda(0, b)\,\,\rmd{d}a + \int \frac{\rmd{d}b}{\lambda(0, b)}\right)
\\&= \int_0^\infty \sin\left((0)u^2 - \frac{b}{u^2}\right) \rmd{d}u = \int_0^\infty \sin\left(-\frac{b}{u^2}\right) \rmd{d}u
\\&= -\sqrt{b}\sqrt{\frac{\pi}{2}} \exp(0)
\end{align*}

$$\implies C_3 =-\sqrt{b}\sqrt{\frac{\pi}{2}},\hspace{0.5cm} \int \lambda(0, b)\,\,\rmd{d}a + \int \frac{\rmd{d}b}{\lambda(0, b)} = 0$$

Let $\lambda(0, b) = f(b)$

$$\int f(b) \,\,\rmd{d}a + \int \frac{\rmd{d}b}{f(b)} = 0$$
$$af(b) + \int \frac{\rmd{d}b}{f(b)} = \textup{constant}$$
$$af(b)= -\int \frac{\rmd{d}b}{f(b)} + \textup{constant}$$

Taking the derivative at $b$,
$$\frac{\rmd{d}}{\rmd{d}b}\left(af(b)\right) = \frac{\rmd{d}}{\rmd{d}b}\left(-\int \frac{\rmd{d}b}{f(b)}\right)$$

$$\implies af'(b) = -\frac{1}{f(b)}$$

$$\implies f'(b)f(b) = \frac{-1}{a}$$
$$\implies f(b)\frac{\rmd{d}\left(f(b)\right)}{\rmd{d}b} = \frac{-1}{a}$$
$$\implies \frac{1}{2}\rmd{d}\left(f(b)^2\right) = \frac{-\rmd{d}b}{a}$$
$$\implies \frac{1}{2} \int \rmd{d}\left(f(b)^2\right)  = \int \frac{-\rmd{d}b}{a}$$

For $\lambda_2, \gamma_2$ constant,
$$\implies \frac{f(b)^2}{2} =  -\frac{b}{a} + \lambda_2$$
$$\implies f(b) = \sqrt{\gamma_2 - \frac{2b}{a}} = \sqrt{\frac{a\gamma_2 - 2b}{a}}$$ \\ \\

Since $f(b) = \lambda(0, b)$ and $f(a) = \lambda(a, 0)$,

$$\implies \lambda(0, b) = \sqrt{\frac{a\gamma_2 - 2b}{a}} = \phi_1\sqrt{\phi_2 - 2b}\,\,\,\textsf{and}\,\,\, \lambda(a, 0) = \sqrt{\frac{b}{b\gamma_1 - 2a}} = \frac{\phi_3}{\sqrt{\phi_4 - 2a}}$$

$\lambda(a, b)$ can be written from $\lambda(0, b)$ and $\lambda(a, 0)$
$$\implies \lambda(a, b) = \phi \sqrt{\frac{\phi_2 - 2b}{\phi_4 - 2a}}\,\,\,\,\,\, \ni \,\, \phi = \frac{\phi_3}{\sqrt{\phi_2}} \,\,\, \textsf{or}\,\,\, \phi = \phi_1\sqrt{\phi_4}.$$

$\phi \neq 0$ since we are looking for non-trivial solutions of the PDE. 

$$\implies \lambda(a, b) = \phi \sqrt{\frac{\phi_2 - 2b}{\phi_4 - 2a}},\hspace{1cm}\phi \neq 0$$

\begin{align*}
F(a, b) &= C_3\exp\left(\phi\int \sqrt{\frac{\phi_2 - 2b}{\phi_4 - 2a}} \rmd{d}a + \int \frac{1}{\phi}\sqrt{\frac{\phi_4 - 2b}{\phi_2 - 2a}} \rmd{d}b\right)
\\&= C_3\exp\left(-\phi\sqrt{\left(\phi_2 - 2b\right)\left(\phi_4 - 2a\right)} - \frac{1}{\phi}\sqrt{\left(\phi_2 - 2b\right)\left(\phi_4 - 2a\right)}\right)
\\&= C_3\exp\left(-\left(\phi + \frac{1}{\phi}\right)\sqrt{\left(\phi_2 - 2b\right)\left(\phi_4 - 2a\right)}\right)
\end{align*}

Since $\phi \neq 0,\,\, \phi + \frac{1}{\phi}  \neq 0$. So, let $\delta = \phi + \frac{1}{\phi} \,\, \ni \,\, \delta \neq 0.$ \\

$$\implies F(a, b) = C_3\exp\left(-\delta\sqrt{\left(\phi_2 - 2b\right)\left(\phi_4 - 2a\right)}\right)$$

\begin{align*}
F(a, 0) &= C_3\exp\left(-\delta\sqrt{\left(\phi_2\right)\left(\phi_4 - 2a\right)}\right)
\\&= \frac{1}{2\sqrt{a}} \sqrt{\frac{\pi}{2}} \exp(0)
\end{align*}

$$\implies C_3 = \frac{1}{2\sqrt{a}},\,\,\,\, \phi_2 = 0.$$
$$\phi_4 \neq 2a,\,\,\, \textsf{since} \,\,\, \phi_4\,\,\, \textsf{is an arbitrary constant, it cannot be a function of}\,\,\, a.$$

\begin{align*}
F(0, b) &= C_3\exp\left(-\delta\sqrt{\left(-2b\right)\left(\phi_4 - 2(0)\right)}\right)
\\&= C_3\exp\left(-\delta\sqrt{\left(-2b\right)\left(\phi_4\right)}\right)
\\&= -\sqrt{b} \sqrt{\frac{\pi}{2}} \exp(0)
\end{align*}

$\implies C_3 = -\sqrt{b},\,\,\, \sqrt{\phi_4 (-2b)} = 0$ $\implies \phi_4 = 0,$ since $b \neq 0$. $b \neq 0$ since $F(a, b) = C_3$ does not satisfy the PDE.

$$\implies F(a,b) =   \frac{1}{2\sqrt{a}} \sqrt{\frac{\pi}{2}}\exp\left(-\delta\sqrt{-2b\left(\phi_4 - 2a\right)}\right)$$
or
$$F(a, b) = -\sqrt{b} \sqrt{\frac{\pi}{2}}\exp\left(-\delta\sqrt{-2a\left(\phi_2 - 2b\right)}\right)$$\\

Determining which of the two satisfies the PDE, let 
$$\displaystyle F(a,b) =   \frac{1}{2\sqrt{a}} \sqrt{\frac{\pi}{2}}\exp\left(-\delta\sqrt{-2b\left(\phi_4 - 2a\right)}\right)$$

\begin{align*}
&\implies \frac{\partial^2}{\partial a \partial b}\left(\frac{1}{2\sqrt{a}} \sqrt{\frac{\pi}{2}}\exp\left(-\delta\sqrt{-2b\left(\phi_4 - 2a\right)}\right)\right) 
\\&\hspace{1.5cm} = \frac{1}{2}\sqrt{\frac{\pi}{2}}\left(\frac{\delta \exp\left(\delta(-\sqrt{4ab - 2b\phi_4})\right)\left(4a\delta\sqrt{b(2a - \phi_4)} - \sqrt{2}\phi_4\right)}{4a\sqrt{a}\sqrt{b(2a - \phi_4)}}\right).
\end{align*}

At $\phi_4 = 0$,
\begin{align*}
\frac{\partial^2}{\partial a \partial b}\left(\frac{1}{2\sqrt{a}} \sqrt{\frac{\pi}{2}}\exp\left(-2\delta\sqrt{ab}\right)\right) &= \frac{1}{2}\sqrt{\frac{\pi}{2}}\left(\frac{\delta \exp\left(-\delta\sqrt{4ab}\right)\left(4a\delta\right)}{4a\sqrt{a}}\right) 
\\&= \frac{1}{2}\sqrt{\frac{\pi}{2}}\left(\frac{\delta^2 \exp\left(-\delta\sqrt{4ab}\right)}{\sqrt{a}}\right) = \frac{\delta^2}{2\sqrt{a}} \sqrt{\frac{\pi}{2}}\exp\left(-2\delta\sqrt{ab}\right)
\end{align*}

If $\delta = 1$, then $\displaystyle\frac{\partial^2}{\partial a \partial b}\left(\frac{1}{2\sqrt{a}} \sqrt{\frac{\pi}{2}}\exp\left(-2\delta\sqrt{ab}\right)\right) =  \frac{1}{2\sqrt{a}} \sqrt{\frac{\pi}{2}}\exp\left(-2\sqrt{ab}\right)$.\\

Also, $\displaystyle F(a, b) =  \frac{1}{2\sqrt{a}} \sqrt{\frac{\pi}{2}}\exp\left(-2\sqrt{ab}\right).$\\

It follows that $\displaystyle F(a, b) =  \frac{1}{2\sqrt{a}} \sqrt{\frac{\pi}{2}}\exp\left(-2\sqrt{ab}\right)$ satisfies the PDE. \\

Also let $\displaystyle  F(a, b) = -\sqrt{b} \sqrt{\frac{\pi}{2}}\exp\left(-\delta\sqrt{-2a\left(\phi_2 - 2b\right)}\right)$
\begin{align*}
&\implies\frac{\partial^2}{\partial a \partial b}\left(-\sqrt{b} \sqrt{\frac{\pi}{2}}\exp\left(-\delta\sqrt{-2a\left(\phi_2 - 2b\right)}\right)\right) 
\\&\hspace{1cm} = \sqrt{\frac{\pi}{2}}\left(\frac{\delta \exp\left(-\delta\sqrt{-2a\left(\phi_2 - 2b\right)}\right)\left(\sqrt{2}\left(4b - \phi_2\right) - 4b\delta\sqrt{a(2b - \phi_2)}\right)}{4\sqrt{b}\sqrt{a\left(2b - \phi_2\right)}}\right).
\end{align*}

At $\phi_2 = 0$,
\begin{align*}
&\frac{\partial^2}{\partial a \partial b}\left(-\sqrt{b} \sqrt{\frac{\pi}{2}}\exp\left(-2\delta\sqrt{ab}\right)\right) = \left(-\delta^2\sqrt{b}\exp\left(-2\delta\sqrt{ab}\right) + \frac{\delta}{\sqrt{a}}\exp\left(-2\delta\sqrt{ab}\right)\right).
\end{align*}

But
\begin{align*}
&\sqrt{\frac{\pi}{2}}\left(-\delta^2\sqrt{b}\exp\left(-2\delta\sqrt{ab}\right) + \frac{\delta}{\sqrt{a}}\exp\left(-2\delta\sqrt{ab}\right)\right)  \neq -\sqrt{b}\sqrt{\frac{\pi}{2}}\exp\left(-2\delta\sqrt{ab}\right)
\end{align*}
for constant values of $\delta$ and this implies that $\displaystyle F(a, b) \neq  -\sqrt{b} \sqrt{\frac{\pi}{2}}\exp\left(-2\delta\sqrt{ab}\right)$ does not satisfy the PDE for all constant values of $\delta$.  \\ \\

Hence,
$$F(a, b) =  \frac{1}{2\sqrt{a}} \sqrt{\frac{\pi}{2}}\exp\left(-2\sqrt{ab}\right)$$\\

By the result in \eqref{speci}, 
$$K_{\frac{1}{2}}\left(x\right)  = \sqrt{\frac{\pi}{2 x}} e^{-x}$$

\begin{align*}
\implies F(a, b) &=  \frac{1}{2\sqrt{a}} \sqrt{\frac{\pi}{2}}\exp\left(-2\sqrt{ab}\right) = \sqrt{\sqrt{\frac{b}{2a}}} \sqrt{\frac{\pi}{2 \cdot 2\sqrt{ab}}} e^{-2\sqrt{ab}} 
\\&= \left(\frac{b}{2a}\right)^{\frac{1}{4}}K_{\frac{1}{2}}\left\{2\sqrt{ab}\right\}
\end{align*}

$$\therefore \int_0^\infty \sin\left(au^2 - \frac{b}{u^2}\right) \rmd{d}u = \left(\frac{b}{2a}\right)^{\frac{1}{4}} K_{\frac{1}{2}}\left\{2\sqrt{ab}\right\}$$
\end{proof} \newpage

\section{Main Results/Propositions}
\begin{tcolorbox}[text width=12cm, lifted shadow={1mm}{-2mm}{3mm}{0.1mm}{black!50!white}]
\begin{align*}
J_{\alpha}(x) &= \frac{1}{\Gamma(\alpha)}\left(\frac{x}{2}\right)^{\alpha} \int_0^1 J_{0}\left(x\sqrt{1 - t}\right)t^{\alpha - 1} \rmd{d}t
\\&=\frac{1}{\pi x\Gamma(\alpha)}\left(\frac{x}{2}\right)^{\alpha} \int_0^{\frac{\pi}{2}} \int_0^\pi 2\sin^{2\alpha - 1}{\var}\sin{\phi}\frac{\partial}{\partial \left(x\cos{\var}\right)}\left(x\cos{\var}\sin\left(x\cos{\var}\sin{\phi}\right)\right) \rmd{d}\phi\rmd{d}\var
\\ I_{\alpha}(x) &= \frac{1}{\Gamma(\alpha)}\left(\frac{x}{2}\right)^{\alpha} \int_0^1 I_{0}\left(x\sqrt{1 - t}\right)t^{\alpha - 1} \rmd{d}t
\\&= \frac{1}{\pi x\Gamma(\alpha)}\left(\frac{x}{2}\right)^{\alpha} \int_0^{\frac{\pi}{2}} \int_0^\pi 2\sin^{2\alpha - 1}{\var}\sin{\phi}\frac{\partial}{\partial \left(x\cos{\var}\right)}\left(x\cos{\var}\sinh\left(x\cos{\var}\sin{\phi}\right)\right) \rmd{d}\phi\rmd{d}\var
\\J_0(x) &= \frac{1}{\pi x}\int_0^{\pi}\sin{\phi}\frac{\partial}{\partial x}\left(x\sin\left(x\sin{\phi}\right)\right)\rmd{d}\phi
\\I_0(x) &= \frac{1}{\pi x}\int_0^{\pi}\sin{\phi}\frac{\partial}{\partial x}\left(x\sinh\left(x\sin{\phi}\right)\right)\rmd{d}\phi
\\&\int_0^\infty \sin\left(au^2 - \frac{b}{u^2}\right) \rmd{d}u = \left(\frac{b}{2a}\right)^{\frac{1}{4}} K_{\frac{1}{2}}\left\{2\sqrt{ab}\right\}
\end{align*}
\end{tcolorbox}

\newpage
\subsection{The area between each of the curves and the $x$-axis gives the values of $J_0(z)$ from $z=1$ to $10$}
\includegraphics[scale=0.68]{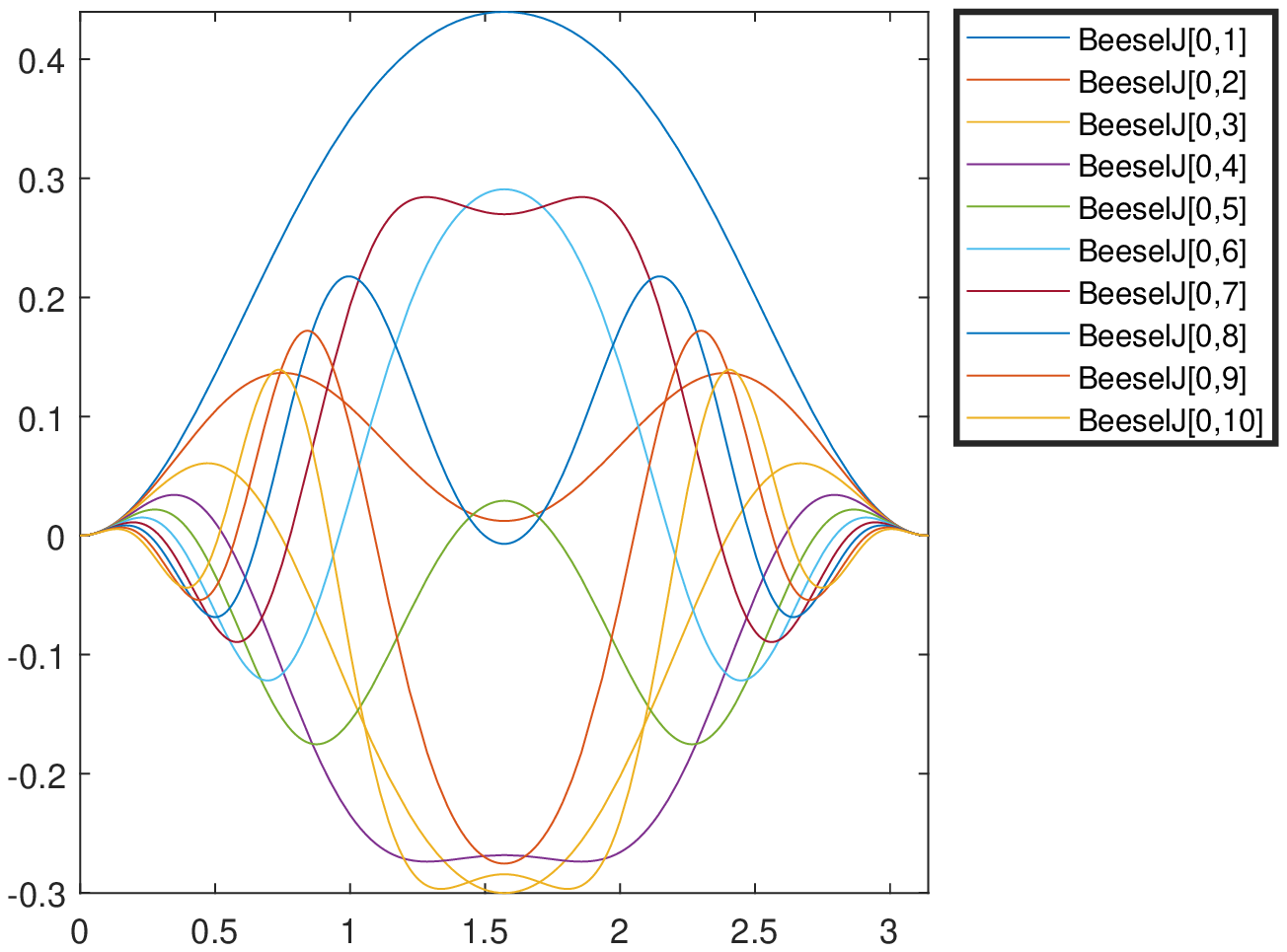}

\subsection{The area between each of the curves and the $x$-axis gives the values of $I_0(z)$ from $z=1$ to $10$}
\includegraphics[scale=0.68]{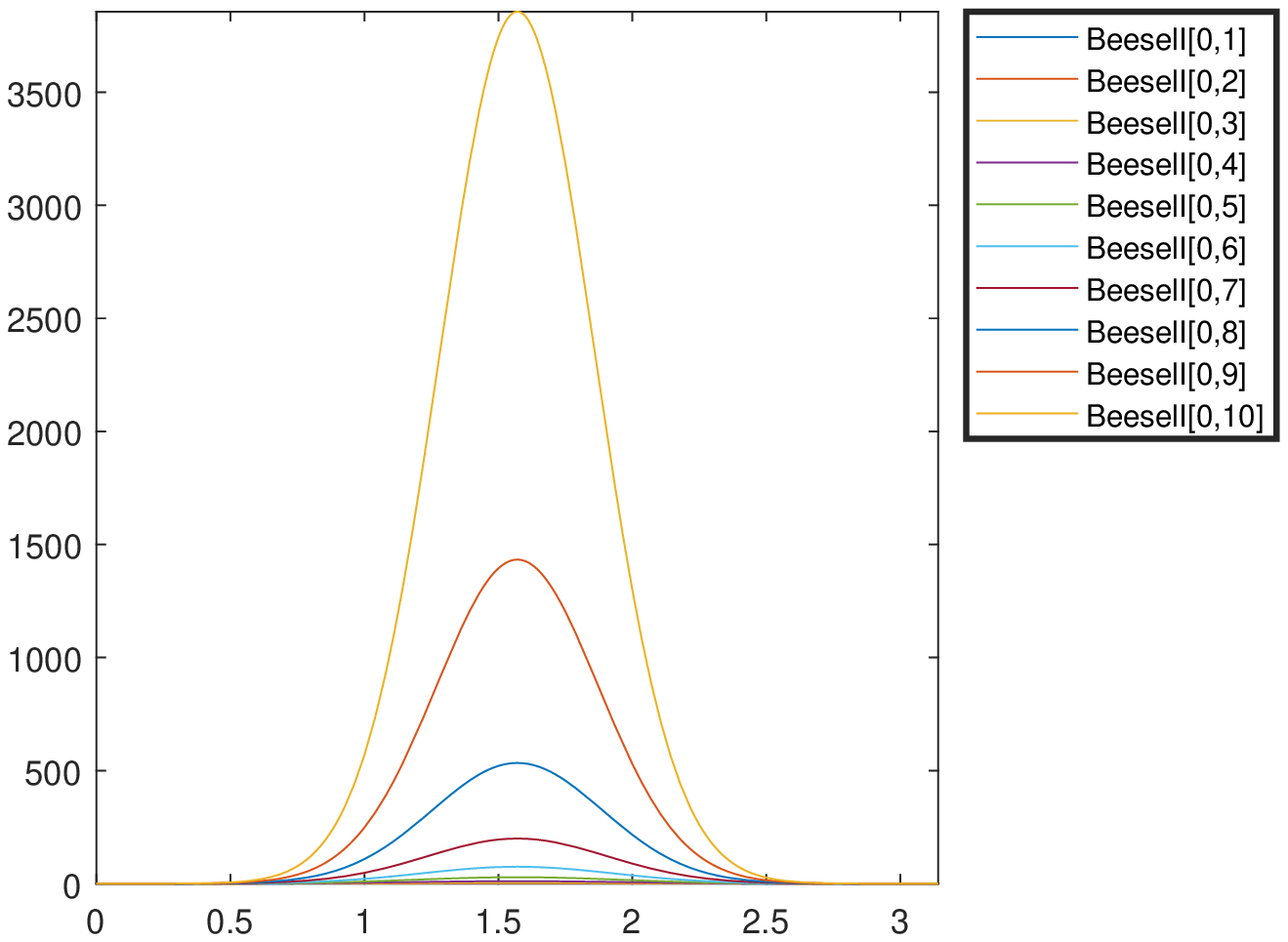}

\newpage
\bibliographystyle{plain}
\bibliography{Bessel}

\begin{thebibliography}{1}

\bibitem{Treat2}
George {N}eville {W}atson~{S}c.{D,} {F}.{R}.{S}.
\newblock {\em A {T}{R}{E}{A}{T}{I}{S}{E} {O}{N} {T}{H}{E} {T}{H}{E}{O}{R}{Y}
  {O}{F} {B}{E}{S}{S}{E}{L} {F}{U}{N}{C}{T}{I}{O}{N}{S}}, chapter~I, page~1.
\newblock Cambridge University Press, 1944.

\bibitem{Treat3}
George {N}eville {W}atson~{S}c.{D,} {F}.{R}.{S}.
\newblock {\em A {T}{R}{E}{A}{T}{I}{S}{E} {O}{N} {T}{H}{E} {T}{H}{E}{O}{R}{Y}
  {O}{F} {B}{E}{S}{S}{E}{L} {F}{U}{N}{C}{T}{I}{O}{N}{S}}, chapter XIII, page
  384.
\newblock Cambridge University Press, 1944.

\bibitem{Treat4}
George {N}eville {W}atson~{S}c.{D,} {F}.{R}.{S}.
\newblock {\em A {T}{R}{E}{A}{T}{I}{S}{E} {O}{N} {T}{H}{E} {T}{H}{E}{O}{R}{Y}
  {O}{F} {B}{E}{S}{S}{E}{L} {F}{U}{N}{C}{T}{I}{O}{N}{S}}, chapter~VI, page 172.
\newblock Cambridge University Press, 1944.

\bibitem{Treat1}
George {N}eville {W}atson~{S}c.{D,} {F}.{R}.{S}.
\newblock {\em A {T}{R}{E}{A}{T}{I}{S}{E} {O}{N} {T}{H}{E} {T}{H}{E}{O}{R}{Y}
  {O}{F} {B}{E}{S}{S}{E}{L} {F}{U}{N}{C}{T}{I}{O}{N}{S}}, chapter~VI, page 180.
\newblock Cambridge University Press, 1944.

\bibitem{Treat5}
\url{https://en.wikipedia.org/wiki/Bessel_function}.

\end{thebibliography}

\end{document}